\documentclass[num-refs]{wiley-article}
\usepackage{graphics}
\usepackage{siunitx}
\usepackage{caption,subcaption}
\usepackage[linesnumbered,ruled,vlined]{algorithm2e}
\makeatletter
\patchcmd{\ttlh@hang}{\parindent\z@}{\parindent\z@\leavevmode}{}{}
\patchcmd{\ttlh@hang}{\noindent}{}{}{}
\makeatother

\papertype{Original Article}

\title{In-situ adaptive reduction of nonlinear multiscale structural dynamics models}

\author[1]{Wanli He}
\author[2]{Philip Avery}
\author[1, 2, 3]{Charbel Farhat}

\affil[1]{Department of Mechanical Engineering, Stanford University, Building 530, 440 Escondido Mall, Stanford, CA 94305-3030, USA}
\affil[2]{Department of Aeronautics and Astronautics, Stanford University, Durand Building, 496 Lomita Mall, Stanford, CA 94305-4035, USA}
\affil[3]{Institute for Computational and Mathematical Engineering, Huang Engineering Center, 475 Via Ortega, Suite 060, Stanford, CA 94305-4042, USA}

\corraddress{Philip Avery, Department of Aeronautics and Astronautics, Stanford University, Durand Building, 496 Lomita Mall, Stanford, CA 94305-4035, USA}
\corremail{pavery@stanford.edu}

\fundinginfo{Air Force Office of Scientific Research, Grant FA9550-17-1-0182; King Abdulaziz City for Science and Technology (KACST)}

\runningauthor{HE \textsc{et al}}

\begin{document}
\maketitle
\begin{abstract}
Conventional offline training of reduced-order bases in a predetermined region of a parameter space leads
to parametric reduced-order models that are vulnerable to extrapolation. This vulnerability manifests
itself whenever a queried parameter point lies in an unexplored region of the parameter space. This paper
addresses this issue by presenting an in-situ, adaptive framework for nonlinear model reduction where
computations are performed by default online, and shifted offline as needed. The framework is based on
the concept of a database of local Reduced-Order Bases (ROBs), where locality is defined in the parameter
space of interest. It achieves accuracy by updating on-the-fly a pre-computed ROB, and approximating
the solution of a dynamical system along its trajectory using a sequence of most-appropriate local ROBs.
It achieves efficiency by managing the dimension of a local ROB, and incorporating hyperreduction in
the process. While sufficiently comprehensive, the framework is described in the context of dynamic multiscale
computations in solid mechanics. In this context, even in a nonparametric setting of the macroscale problem
and when all offline, online, and adaptation overhead costs are accounted for, the proposed computational
framework can accelerate a single three-dimensional, nonlinear, multiscale computation by an order of
magnitude, without compromising accuracy.
\keywords{dynamics, hyperreduction, local reduced-order basis, multiscale, nonlinear model reduction, parameter sampling}
\end{abstract}

\section{Introduction}

Heterogeneous materials are encountered to an increasing extent in simulation-based engineering applications
involving structural analysis of solids. Applications include, but are not limited to, advanced modular
armor protection \cite{qiao2008impact}, design of new materials with special-purpose optimized properties
\cite{ma2017towards}, and bone fracture risk assessment \cite{hoffler2000heterogeneity}. When the scale
of the heterogeneity is such that multiple disparate spatial scales require resolution in order to achieve
a reliable prediction, multiscale homogenization can be a viable alternative to prohibitively expensive
direct numerical simulation involving a massive monolithic discretization \cite{mei2010homogenization}. 

In principle, the method of multiscale homogenization assumes that the material configuration is homogeneous
at the macroscale level, but heterogeneous at the level of the smallest represented scale. The evaluation
of a conventional constitutive response function at each macroscale material point is replaced by the
solution of a discrete Boundary Value Problem (BVP) whose computational domain is a so-called microscale
Representative Volume Element (RVE) characterizing the microstructure. The complete BVP is formulated
with prescribed boundary deformations corresponding to the deformation gradient at the macroscale material
point. The sequence of scales is coupled by formalized procedures of localization and homogenization
to transmit information between differing scales \cite{smit1998prediction, miehe1999computational,feyel2000fe2,
kouznetsova2001approach}. In this work, a particular variant of the method of multiscale homogenization
known as ``Finite Element squared'' or FE$^2$ is considered as a starting point. As the name implies,
in this variant the Finite Element (FE) method is used to discretize both the macroscopic computational
domain and the microscale RVE. The repeated solution of the FE BVP at each material point is an acknowledged
computational bottleneck, motivating the development of alternative surrogate models to accelerate this
task. Projection-based Model Order Reduction (PMOR) equipped with a second-tier hyperreduction approximation
is one established technique for reducing the computational complexity of nonlinear FE models that has
recently been demonstrated to be an effective candidate surrogate model in the context of multiscale
homogenization \cite{zahr2017multilevel}. Alternative surrogate model candidates include regression-based
models constructed using various methods including Kriging and artificial neural networks. Such methods
can potentially lead to more computationally efficient surrogate model evaluation compared to PMOR; however
a large amount of data may be required to train the model when the constitutive response is complex.
To accelerate the acquisition of training data, PMOR again presents a useful computational technology.
In this work, a multifidelity framework leveraging low-, medium- and high-fidelity models is proposed to
achieve a computationally tractable end-to-end procedure for nonlinear multiscale simulations. The proposed
framework features a novel in-situ construction and utilization of an adaptive PMOR to eliminate computation
bottlenecks. The usefulness of this novelty is not restricted to the current application of interest
(i.e. multiscale modeling); benefits are anticipated in other related processes such as multifidelity
optimization. It is therefore presented in what follows as a general-purpose methodology.

PMOR is typically embedded in an offline-online computational framework. In the phase of non-time-critical,
offline training, information data from a High Dimensional Model (HDM) are collected and used to construct
a Reduced-Order Basis (ROB) -- for example, using Proper Orthogonal Decomposition (POD) and the method
of snapshots \cite{sirovich1987turbulence} -- and the corresponding Projection-based Reduced Order Model
(PROM). Then, in the phase of time-critical, online evaluation, an approximate solution is computed using
the PROM as a surrogate model for the underlying HDM. This offline-online decomposition of the computational
effort is appropriate in real-time applications that are performed for the exploration of a space of
interest associated with some usually but not necessarily physical model parameters, such as in design
optimization \cite{legresley2000airfoil, manzoni2012shape}, uncertainty quantification \cite{farhat2018stochastic,
farhat2018modeling}, and repeated analysis \cite{hetmaniuk2012review,amsallem2010towards}. In these cases
where the space of one or more physical macroscopic parameters (such as those defining the geometry, boundary/initial
conditions, or material properties) is explored, the offline cost is amortized by the sizable number
of required analyses associated with the number of parameter points necessary to adequately perform this
exploration. Consequently, the overall speedup can still be significant even when the potentially substantial
offline cost is accounted for. Here, it is emphasized that in the context of multiscale homogenization,
PMOR is computationally feasible even in the nonparametric setting due to the large number of microscale
RVE simulations that are required at a single macroscopic parameter point. Hence, a desirable speedup
of the multiscale computations can be still achieved by building PROMs for all smaller scales. Furthermore,
a nonparametric application calls for an \emph{in-situ} alternative to the conventional offline-online
strategy, where the construction and utilization of the hyperreduced PROM (HPROM) can be performed during
the multiscale simulation without requiring any preliminary computational work, using solution snapshots
collected during already performed microscale RVE computations.

Several existing PMOR methods have been previously developed and applied to the solution of nonlinear
multiscale problems. In \cite{yvonnet2007reduced}, a PROM of the microscale problem was constructed based
on an \'a priori training strategy that involved computing solution snapshots using predefined microstrain
histories. This approach was subsequently extended in \cite{monteiro2008computational} to formulations
with periodic boundary conditions for the microscale model. However, in each of these works the PROMs
were not equipped with a hyperreduction method and consequently the resulting speedup factor was limited
to about an order of magnitude, even without accounting for the offline costs. More recently, the co-authors
of \cite{zahr2017multilevel} introduced a nonlinear PMOR framework incorporating the Energy Conserving
Sampling and Weighting (ECSW) hyperreduction method \cite{farhat2014dimensional,farhat2015structure}
which significantly reduced the computational cost of multiscale problems characterized by large deformations
and material nonlinearity, including when the offline costs were taking into account. Specifically, the
computational framework proposed in \cite{zahr2017multilevel} included a two-step training strategy where:
a microscale HPROM is first constructed in-situ in order to achieve significant speedups even in nonparametric
settings; next, a conventional offline-online training approach is performed to build a parametric macroscale
HPROM. A notable feature of this approach is the reduction in the cost of the training by using the in-situ
microscale HPROM to accelerate macroscale snapshot acquisition. However, a weakness of the proposed framework
is its vulnerability to extrapolation if the underlying approximation subspace is not sufficiently explored
during the in-situ training phase.

The alternative in-situ training strategy proposed in this paper avoids extrapolation by incorporating
an adaptive HPROM construction based on the concept of a database of local Reduced-Order Bases (ROBs)
constructed and updated on-the-fly. This process treats the components of the deformation gradient as
a parameter vector, which facilitates the locating of each microscale BVP within the database in order
to assign to it online the most-appropriate local ROB. Accuracy is optimized by collecting new snapshots
and other additional information as needed and updating accordingly the local ROBs. This usage of local
ROBs enables the efficient solution of difficult problems for which PMOR is not effectively performed
by constructing a single global ROB. Note that locality in this paper refers to the parameter space characterized
by the deformation gradient, in sharp contrast with previous works where locality was considered with
respect to the solution manifold \cite{amsallem2012nonlinear}.

The remainder of this paper is organized as follows. Section \ref{sec:theory} provides a brief summary
of the FE$^2$ multiscale homogenization method. Section \ref{sec:pmor} introduces the microscale HPROM
and its training strategy. This section also highlights the necessity for ROB adaptation to avoid extrapolation. 
Section \ref{sec:adaptive-ROM} proposes a novel, online, in-situ ROB adaptation approach suitable for
nonlinear, multifidelity, multiscale, computational homogenization. Section \ref{sec:apps} provides an
assessment of the overall proposed adaptive PMOR framework using several representative applications,
and demonstrates its ability to deliver speedups while compromising neither stability nor accuracy.

\section{The FE$^2$ computational framework}
\label{sec:theory}

In this section a brief summary of the FE$^2$ method is provided, focusing on the aspects that are relevant
to subsequent developments including notation. For a comprehensive formulation of multiscale homogenization,
the interested reader is referred to \cite{zahr2017multilevel} and the references cited therein. The
fundamental generalization underlying FE$^2$ is that the evaluation of stress as a function of strain
for a heterogeneous material is governed by the solution of a locally attached boundary value problem
characterizing its microstructure. This concept can be applied recursively to materials that exhibit
$n_s+1$ separated scales; in what follows $k=0, \ldots, n_s$ denotes the sequence of scales (or levels),
with $k=0$ and $k=1, \ldots, n_s$ designating the macroscale and the finer meso/microscales, respectively.
Specifically, at the $n_s$ coarsest scales ($k=0, \ldots, n_s-1$) each stress-strain relation is associated
with a locally attached microstructure, while at the finest scale ($k=n_s$) the stress-strain relation
is defined by a conventional constitutive law.

\subsection{Semi-discrete governing equations}
\label{sec:governing}

In this work, a distinction is made between constrained and unconstrained Degrees Of Freedom (DOFs),
the former being those whose values are specified by essential boundary conditions, while the latter
are subject to no such restriction. Following the notation of \cite{zahr2017multilevel}, $\stackrel{\circ}{\boldsymbol{v}}_k$
is defined here as a vector over the set of constrained DOFs at level $k$, $\boldsymbol{v}_k$ as the
corresponding vector over the set of unconstrained DOFs at level $k$, and $\overline{\boldsymbol{v}}_k$
as the concatenation of these two vectors. A similar notation is adopted for matrices. It follows that
a subscript $0$ attached to any vector or matrix designates a macroscopic quantity. Using this notation
and assuming that $\ddot{\stackrel{\circ}{\boldsymbol{u}}}_{0}=0$, the semi-discrete form of the equations
of motion governing dynamic equilibrium can be written as
\begin{equation}
\label{eq:governing-semi}
\begin{aligned}
  \boldsymbol{M}_{0}\ddot{\boldsymbol{u}}_{0}+\boldsymbol{f}_{0}\left(\boldsymbol{u}_{0},
  \stackrel{\circ}{\boldsymbol{u}}_{0}, \boldsymbol{\xi}_{0}\right) &= \boldsymbol{0} \\ 
  \boldsymbol{f}_{k}\left(\boldsymbol{u}_{k}, \stackrel{\circ}{\boldsymbol{u}}_{k}, \boldsymbol{\xi}_{k}\right)
  &= \boldsymbol{0}
\end{aligned}
\end{equation}
for $k=1,\ldots,n_s$, where $\boldsymbol{M}$ is a FE mass matrix, $\boldsymbol{f}$ is a FE vector of
internal forces, $\boldsymbol{u}$ and $\ddot{\boldsymbol{u}}$ are FE vectors of displacements and accelerations,
respectively, and $\boldsymbol{\xi}$ is a vector of history variables associated with certain constitutive
laws such as elastoplasticity when introduced at the level $k=n_s$. For levels $0\le k<n_s$, $\boldsymbol{\xi}_{k}$
consists of the set of history variables at all material points of the attached microstructure. The governing
equations \eqref{eq:governing-semi} are coupled by the localization and homogenization conditions presented
in the next section.

\subsection{Scale transmission conditions}
\label{sec:transmission}

Kinematic information is transmitted from coarser to finer scales via prescribed nonhomogeneous displacement
boundary conditions applied to the exterior boundary of the attached microstructure. These prescribed
values are specified such that the volumetric average of the deformation gradient at the finer scale
is identical to the deformation gradient at the material point of the coarser scale to which the microstructure
is attached. This is referred to as the localization condition and can be expressed as follows
\begin{equation}
\label{eq:localization}
  \stackrel{\circ}{\boldsymbol{U}}_{k+1}=\stackrel{\circ}{\boldsymbol{X}}_{k+1}\left(\boldsymbol{F}_k
  -\boldsymbol{I}\right)
\end{equation}
where $\boldsymbol{F}_k$ is the coarse-scale deformation gradient, $\boldsymbol{I}$ is the $3 \times 3$
identity matrix, and $\stackrel{\circ}{\boldsymbol{U}}_{k+1}$ is the matrix whose three columns contain
the prescribed $x$-, $y$- and $z$-components of displacement. In other words, the displacement vector
over the set of constrained DOFs at level $k+1$, $\stackrel{\circ}{\boldsymbol{u}}_{k+1}$, is obtained
-- up to a permutation -- by stacking the columns of $\stackrel{\circ}{\boldsymbol{U}}_{k+1}$. Similarly,
$\stackrel{\circ}{\boldsymbol{X}}_{k+1}$ is the matrix whose three columns contain the $x$-, $y$- and
$z$-coordinates in the reference configuration of the nodes on the exterior boundary of the attached
microstructure whose DOFs are constrained. In the particular case of the so-called uniform essential
boundary conditions, the values of all DOFs located on the exterior boundary of the attached microstructure
are prescribed. Alternatively, in the case of periodic boundary conditions, only a minimal set of corner
points have prescribed displacements while the remaining exterior boundary DOFs are constrained using
linear multipoint constraint equations to enforce the periodicity of the fluctuations between each pair
of opposite faces.

After computing the response of an attached microstructure to the prescribed boundary displacements,
information is transmitted back from finer to coarser scales in the form of a homogenized stress tensor
defined as the volumetric average of the stress in the attached microstructure. This quantity can be
conveniently obtained from the reaction forces associated with the constrained boundary DOFs as follows
\begin{equation*}
\label{eq:homogenization}
  \boldsymbol{P}_{k}=\frac{1}{\left|\mathcal{B}_{k+1}\right|} \boldsymbol{X}_{k+1}^{T} \stackrel{\circ}{\boldsymbol{R}}_{k+1}
\end{equation*}
where $\boldsymbol{P}_{k}$ is the coarse-scale homogenized first Piola-Kirchhoff stress tensor, $\left|\mathcal{B}_{k+1}\right|$
is the volume of the RVE, and $\stackrel{\circ}{\boldsymbol{R}}_{k+1}$ is the matrix whose three columns
contain the $x$-, $y$-, and $z$-components of the internal force vector of the nodes on the exterior
boundary of the attached microstructure whose DOFs are constrained, i.e. the boundary reaction forces.

In summary, the multiscale computational model distinguishes itself from a conventional FE model in that
it involves repeatedly solving one or more microscale BVPs, with one such solution required for each
stress-strain evaluation at all scales other than the finest. For a given RVE, the microscale BVPs vary
only in the values of the prescribed boundary displacements. Regardless of whether uniform or periodic
boundary conditions are adopted, \eqref{eq:localization} shows that the computational response of the
microscale model can be characterized as an implicit function of nine components of the deformation gradient
$\boldsymbol{F}_{k-1}$. In Section \ref{sec:adaptive-ROM}, the set of components of $\boldsymbol{F}_{k-1}$
is treated as a parameter vector as in the parametric setting of PMOR. The parameter space $\mathcal{D}_{k}
\subset \mathbb{R}^{n_{sd}}$, where $n_{sd} = 4$ for two-dimensional problems, $n_{sd} = 9$ for three-dimensional
ones, and $\boldsymbol{F}_{k-1} \in \mathcal{D}_{k}$, is explored by solving the microscale BVP repeatedly
in order to generate snapshots of the microstructural response as required for the construction of a
surrogate model.

At each computational time-step of a single multiscale simulation, a total of $\prod_{k=0}^{n_s-1}\left|\mathcal{G}_{k}\right|$ 
high-dimensional, microscale BVPs are solved at the finest scale, where $\left|\mathcal{G}_{k}\right|$
denotes the number of Gauss points at level $k$. The overall computational cost associated with such
a computational complexity can be prohibitive for large-scale multiscale models, which motivates PMOR
to truncate the computational complexity by reducing the operational dimensionality of the multiscale
model. Even in the absence of a physical parametric setting at the macroscale level, computational savings
are still feasible if the macroscale model is not reduced, but rather a sequence of HPROMs are constructed
at all other scales $k$, $k=1,\ldots,n_s$. For this setting, the co-authors of \cite{zahr2017multilevel}
proposed an in-situ training strategy for microscale reduction where only the single multiscale simulation
of interest is required to be performed, without any predefined offline computational work. This methodology
will be described in the following section, after which a novel contribution will be presented to extend
the scope and robustness of the in-situ concept without sacrificing its inherent convenience.

\section{Projection-based model order reduction}
\label{sec:pmor}

Here, the application of PMOR to the FE$^2$ computation framework is reviewed. Specifically, the POD
reduction and ECSW hyperreduction methods are used to construct a surrogate model to accelerate the solution
of each locally attached microstructural FE model at all fine scales ($k=1, \ldots, n_s$). Training is
performed in-situ (or on-the-fly), pursuant to our objective of enabling a computational technology for
multiscale modeling that is tractable, robust, and convenient for simulators accustomed to conventional
materials. Although not considered in what follows, PMOR can also be applied to the macroscale BVP ($k=0$,
see \cite{zahr2017multilevel}).

\subsection{Microscale reduction framework}

At each fine scale $k$, $k=1,\ldots,n_s$, the number of DOFs $n_{k}^d$ of the computational FE model
is reduced by searching for the solution $\boldsymbol{u}_{k}$ of the problem of interest in a low-dimensional
subspace -- that is,
\begin{equation}
\label{eq:basis}
  \boldsymbol{u}_{k} \approx \boldsymbol{V}_{k} \boldsymbol{y}_{k}
\end{equation}
where $\boldsymbol{V}_{k} \in \mathbb{R}^{n_{k}^d \times s_{k}}$ is a ROB constructed using POD and the
method of snapshots \cite{sirovich1987turbulence}, $\boldsymbol{y}_{k} \in \mathbb{R}^{s_{k}}$ is a vector
of generalized coordinates, and $s_{k} \ll n_{k}^d$. The snapshots from which the ROB is constructed
correspond to high-dimensional solutions of the microstructural FE model that are collected and processed
\emph{in-situ}. The term in-situ indicates that such solutions are obtained during the same multiscale
simulation in which the model reduction framework is utilized, in contrast to online-offline treatments.

The dimensionality of the semi-discrete equations governing the microscale BVP \eqref{eq:governing-semi}
is reduced by performing a Galerkin projection -- that is, substituting \eqref{eq:basis} into these equations
and pre-multiplying them by the transpose of $\boldsymbol{V}_{k}$. This leads to the microscale PROM
\begin{equation}
\label{eq:governing-ROM}
  \boldsymbol{V}_{k}^{T}\boldsymbol{f}_{k}\left(\boldsymbol{V}_{k}\boldsymbol{y}_{k}, \stackrel{\circ}{\boldsymbol{u}}_{k},
  \boldsymbol{\xi}_{k}\right) = \boldsymbol{0}, \qquad k=1,\ldots,n_s
\end{equation}
where the superscript $T$ designates the transpose operation.

Furthermore, a second-tier hyperreduction approximation is also introduced to avoid computations that
scale with the size $n_{k}^d$ of the HDM at each level $k > 0$. In this work, the ECSW method \cite{farhat2014dimensional,
farhat2015structure} is preferred for this purpose due to its structure-preserving property \cite{farhat2015structure}
and the simplicity by which its approximation error can be controlled, namely, via a user-specified thresholding
parameter. As the name suggests, ECSW samples a subset $\widetilde{\cal E}_{k}$ of the element set ${\cal E}_{k}$
of a given FE mesh -- that is, $\widetilde{\cal E}_{k} \subset {\cal E}_{k}$ -- and assigns to each sampled
element $e$ a positive weight $\alpha_{k}^{e}\geqslant 0$. The internal force term in the microscale
PROM \eqref{eq:governing-ROM} is evaluated by simply weighting the contributions of the sampled elements
only, i.e.
\begin{equation}
\label{eq:ECSW}
  \boldsymbol{V}_{k}^{T} \boldsymbol{f}_{k}\left(\boldsymbol{V}_{k} \boldsymbol{y}_{k}, \stackrel{\circ}{\boldsymbol{u}}_{k},
  \boldsymbol{\xi}_{k}\right)
	\approx \sum_{e \, \in \, {\widetilde{\cal E}}_{k}} \alpha_{k}^{e}\left(\boldsymbol{V}_{k}^{e}\right)^{T}
  \boldsymbol{f}_{k}^{e}\left(\boldsymbol{V}_{k}^{e} \boldsymbol{y}_{k}, \stackrel{\circ}{\boldsymbol{u}}_{k}^{e},
  \boldsymbol{\xi}_{k}^{e}\right)
  = \boldsymbol{f}_{k}^{r}\left(\boldsymbol{y}_{k}, \stackrel{\circ}{\boldsymbol{u}}_{k}, \boldsymbol{\xi}_{k}\right) 
\end{equation}
In the above expressions, the superscript $e$ designates the restriction of a global vector or matrix
to element $e$. The optimal {\it reduced mesh} $\widetilde{\cal{E}}_{k}$ and associated weights can be
estimated using either the Non-Negative Least Squares (NNLS) algorithm proposed by Lawson and Hanson
\cite{lawson1995solving} or alternative sparsity promoting L1-minimization methods \cite{chapman2017accelerated}.

In the context of the microscale reduction framework, reduction and hyperreduction are considered not
only for the unconstrained part of the internal force vector as shown in \eqref{eq:ECSW}, but also for
both the constrained part of the internal force required for the homogenization condition and the linear
constraint equations required for the enforcement of periodic boundary conditions \cite{zahr2017multilevel}.

\subsection{In-situ training strategy}
\label{sec:rom-micro}

As already stated in Section \ref{sec:theory}, the deformation gradient $\boldsymbol{F}_{k-1}$ at a Gauss
point characterizes the microscale BVP at level $k$, $k=1,\ldots,n_s$. In what follows, each microscale
BVP at level $k$ is associated with a parameter vector $\boldsymbol{z}_{k}$ in the parameter space $\mathcal{D}_{k}
\subset \mathbb{R}^{sd}$ (also referred to in the remainder of this section as parameter ``point'' in
$\mathcal{D}_{k}$), where $\boldsymbol{z}_{k}$ is obtained by stacking the columns of $\boldsymbol{F}_{k-1}$.
This parametric interpretation of the attached microstructure BVPs motivates the adaptive reduction process
discussed in Section \ref{sec:adaptive-ROM}.

At each computational time-step of a multiscale simulation, $\prod_{j=0}^{k-1}\left|\mathcal{G}_{j}\right|$
microscale BVPs are solved at level $k$, $k=1, \ldots, n_s$. This suggests that sufficient training data
can be acquired for the purpose of constructing a PROM at each level $k$, $k = 1,\ldots,n_s$. An in-situ
microscale training strategy proposed in \cite{zahr2017multilevel} that leverages this observation can
be described by the following steps:
\begin{enumerate} 
  \item Solve at level $k$, $k = 1,\ldots,n_s$, a sufficient number $n_{k-1}^{\mathcal{G}}$ of high-dimensional
        microscale BVPs encountered at the Gauss points of level $k-1$, and sample the resulting solution
        snapshots (note that $n_{k}^{\mathcal{G}}$ is not necessarily equal to $\left|\mathcal{G}_{k}\right|$).  
  \item Construct a global ROB $\boldsymbol{V}_{k}$.
  \item Construct a reduced mesh $\widetilde{\cal{E}}_{k}$.
  \item Exploit the resulting HPROM at level $k$ to solve all remaining microscale BVPs at this level.
\end{enumerate}
This nonadaptive, in-situ, microscale training strategy, where the choice for $n_{k-1}^{\mathcal{G}}$ is
user-specified (and therefore arbitrary), may encounter scenarios involving extrapolation when the remaining
microscale BVPs visit an unexplored region of the parameter space $\mathcal{D}_{k}$. If a queried parameter
point falls outside of the subregion defined by the set of sampled points, a good approximate solution
of the microscale BVP may not be computable in the existing low-dimensional subspace $\boldsymbol{V}_{k}$.
However, an \emph {adaptive} in-situ training strategy has the potential to avoid extrapolation by collecting
new high-dimensional solution snapshots and using them to enrich the subspace of approximation when a
scenario that would otherwise involve extrapolation is detected.

Typically, the enrichment of the approximation subspace needed when the explored domain within $\mathcal{D}_{k}$
is enlarged requires the expansion of the dimension of the ROB, leading to a deterioration in the attainable
speedup. When the dimension of the expanded ROB exceeds a certain threshold, a method is proposed next
to split the global ROB into two or more local ROBs, and subsequently find an approximate solution at
each queried parameter point using the most suitable of pre-computed and stored local ROBs. A similar
approach was proposed in \cite{amsallem2012nonlinear} for constructing local ROBs, where locality was
defined however with respect to the solution manifold. This approach is not ideally-suited for use with
in-situ microscale training as it requires the clustering of solution snapshots, which can be computationally
expensive in this context. Thus, this paper proposes an adaptive PMOR strategy with a method of databases
containing local ROBs, where locality is defined with respect to the parameter space. The details of
this proposal are presented in Section \ref{sec:adaptive-ROM}.

\section{Adaptive in-situ sampling and local reduced-order bases}
\label{sec:adaptive-ROM}

Next, an in-situ, adaptive sampling process based on the concept of a database $\mathcal{DB}_k$ of local
ROBs at each level $k$, $k=1,\ldots,n_s$, is presented. It leverages the accuracy and computational efficiency
of the PMOR framework. The detection of extrapolation when solving a queried microscale BVP is performed
using a residual-based error indicator. When an update criterion is triggered by this indicator at level
$k$, new information is sampled and used to enrich the content of $\mathcal{DB}_k$. In order to preserve
the computational efficiency of PMOR, a local ROB is split into two smaller ones if its current size
exceeds a certain limit. 

In principle, each time a local ROB is updated, the corresponding reduced mesh should also be updated
-- because the training of ECSW depends on the state of the ROB. Therefore, the HPROM should be updated
too accordingly. However, efficiently performing such updates is outside the scope of this paper: it
is the subject of ongoing research. With multiple local ROBs stored in a database $\mathcal{DB}_k$, the
most-appropriate local ROB and associated HPROM must be identified to solve a microscale BVP queried
at level $k$. 

In order to construct such an adaptive PMOR framework, both the microscale solution snapshots and the
underlying deformation gradient parameter points are collected. In order to describe this adaptive framework,
the following additional nomenclature and notation are introduced:
\begin{itemize}
  \item ${\cal S}_z^{m_k} = \left\{\boldsymbol{z}_{k}^{(1)}, \ldots, \boldsymbol{z}_{k}^{\left(m_{k}\right)}\right\}$,
        denotes a set of $m_{k}$ deformation gradient parameter points considered among those that have
        already been sampled at level $k$, and $\boldsymbol{z}_{k}^{(i)} \in \mathcal{D}_{k}$, $i=1,\ldots,m_{k}$.
        It is collected in the matrix
        \begin{equation*}
          \boldsymbol{N}_{k}^{z}=\left[\boldsymbol{z}_{k}^{(1)} \cdots \boldsymbol{z}_{k}^{\left(m_{k}\right)}\right]
          \in \mathbb{R}^{{n_{sd}}\times m_k}
        \end{equation*}
  \item $\boldsymbol{N}_{k}^{u}$ denotes the corresponding snapshot matrix whose $i$-th column stores
        the HDM-based solution of the microscale BVP characterized by $\boldsymbol{z}_{k}^{(i)}$.
\end{itemize}

\subsection{Detection of extrapolation and other large errors}

\begin{figure}[h!] 
  \centering
  \includegraphics[width=0.5\textwidth]{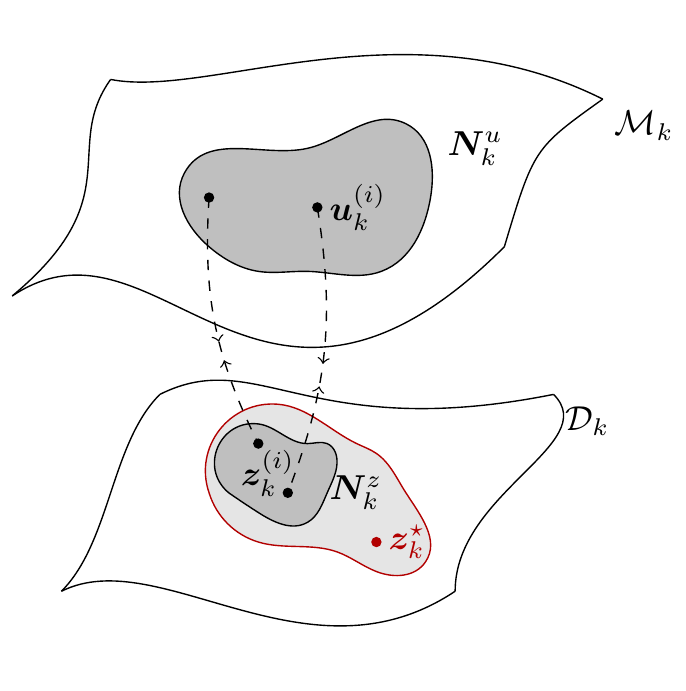} 
  \caption{Extrapolation scenario for an in-situ-trained microscale HPROM: dark/light shaded area designates
           the region of the parameter space explored/unexplored during training.}
  \label{fig:extrap}
\end{figure}

When applied to the solution of a microscale BVP queried at an unsampled parameter point $\boldsymbol{z}_{k}^{\star}
\in \mathcal{D}_{k}$, a local, in-situ-trained, microscale HPROM will perform the equivalent of an extrapolated
approximation if the sampling process underlying the training of the corresponding local ROB did not
explore a region of $\mathcal{D}_{k}$ containing $\boldsymbol{z}_{k}^{\star}$. This scenario is illustrated
in Figure \ref{fig:extrap}, where ${\cal M}_k$ denotes the solution manifold at level $k$. Hence, if
its pitfalls are to be avoided, the potential for this scenario must be detected before the solution
of the queried microscale BVP using a most-appropriate local HPROM is computed or accepted.

There are two main approaches for detecting an extrapolation associated with the approximate solution
of a queried microscale BVP characterized by $\boldsymbol{z}_{k}^{\star}$: a geometric approach that
determines the position of $\boldsymbol{z}_{k}^{\star}$ in $\mathcal{D}_k$, relative to the sampled parameter
points collected in $\boldsymbol{N}_{k}^{z}$; and a residual-based approach that computes the residual
obtained when approximating the HDM-based solution of the queried microscale problem using a most-appropriate
local HPROM.

The geometric approach is conceptually simple. However, its implementation is not trivial. It can be
carried out by computing the hypersphere of smallest radius containing all parameter points collected
in $\boldsymbol{N}_{k}^{z}$, or the minimal hyperpolygon connecting these parameter points, then checking
whether $\boldsymbol{z}_{k}^{\star}$ lies or not within the hypersphere or hyperpolygon. Either implementation
necessitates a fair amount of efficient computational geometry. In either case however, the main drawback
of this approach is that if $\boldsymbol{z}_{k}^{\star}$ falls within either of the aforementioned target
domains, this does not necessarily imply that the approximate solution computed using the most-appropriate
local HPROM has an acceptable accuracy. Such information is provided however by the residual-based approach
which, unlike the geometric approach, checks for potential extrapolation after rather than before the
solution of the queried microscale BVP is computed. Specifically, this approach computes the following
relative residual magnitude 
\begin{equation}
\label{eq:residual}
  r_{k}=\frac{\left\|\boldsymbol{f}_{k}\left(\boldsymbol{V}_{k}\boldsymbol{y}_{k}, \stackrel{\circ}{\boldsymbol{u}}_{k},
  \boldsymbol{\xi}_{k} \right)\right\|}{\left\|\boldsymbol{f}_{k}\left(\boldsymbol{0}, \stackrel{\circ}{\boldsymbol{u}}_{k},
  \boldsymbol{\xi}_{k} \right)\right\|}
\end{equation}
and adopts it as an error indicator. When it exceeds a prescribed tolerance $r_{\text{tol}}$, it suggests
discarding the approximate solution $\boldsymbol{V}_{k}\boldsymbol{y}_{k}$, enriching the most-appropriate
local ROB $\boldsymbol{V}_{k}$, and repeating this process until $r_{k} \leq r_{\text{tol}}$. Hence, this
alternative approach, which leads to the updating criterion outlined in Algorithm \ref{alg:criterion},
equally applies for assessing extrapolation and interpolation errors.

\begin{algorithm}[H]
\DontPrintSemicolon
  \KwData{Queried parameter point $\boldsymbol{z}_{k}^{\star} \in \mathcal{D}_k$, and user-specified
          tolerance value $r_{\text{tol}} > 0$}
  \KwResult{Updating flag $\textrm{UF}$}
  Compute the approximate solution of the microscale BVP characterized by $\boldsymbol{z}_{k}^{\star}$
  using the most-appropriate local HPROM\;
  Compute the corresponding relative residual magnitude $r_{k}^{\star}$ \eqref{eq:residual}\; 
  \eIf{$r_{k}^{\star} > r_{\text{tol}}$}{ 
    Solve the microscale BVP using the corresponding HDM\; 
    Collect the computed solution $\boldsymbol{u}_{k}^{\star}$ in $\boldsymbol{N}_{k}^{u}$\; 
    Sample the parameter point $\boldsymbol{z}_{k}^{\star}$ and collect it in $\boldsymbol{N}_{k}^{z}$\; 
    Set $\textrm{UF} = \textrm{true}$
  }{
    Accept the HPROM solution $\boldsymbol{u}_{k}^{\star}$\; 
    Set $\textrm{UF} = \textrm{false}$
  }
  \caption{Updating criterion for a database of local HPROMs.}
  \label{alg:criterion}
\end{algorithm}

\subsection{Updating of a local reduced-order basis}

Each update of a database $\mathcal{DB}_k$ starts by solving a new microscale BVP using the corresponding
HDM, appending the computed solution $\boldsymbol{u}_{k}^{\star}$ to the current content of the relevant
snapshot matrix $\boldsymbol{N}_{k}^{u}$, performing the Singular Value Decomposition (SVD) of $\boldsymbol{N}_{k}^{u}$,
and deducing a new ROB using the singular value energy truncation criterion. Because the computational
complexity of SVD applied to a rectangular matrix scales quadratically with the number of columns of
the matrix, a low-rank SVD update method such as that presented in \cite{zahr2015progressive} is preferred
here for compressing the collected snapshots.

For a highly nonlinear problem, the SVD and singular value energy truncation criterion configured for
accuracy typically lead to a large ROB. Hence, if the size of such a ROB exceeds a specified capacity
limit $c_{\text{max}}$, this ROB is split here into two smaller ROBs, each of which capable of capturing
only the local behavior of a nonlinear microstructure. This splitting, which is bound in this context
to occur recursively, leads to the concept of a local ROB. This concept is similar to that of local ROBs
pioneered in \cite{amsallem2012nonlinear} and equipped with the notion of clustering of solution snapshots.
However, it differs from it in that locality refers here to the parameter space $\mathcal{D}_k$, whereas
locality in \cite{amsallem2012nonlinear} referred to the solution manifold (denoted here by ${\cal M}_k)$. 
Indeed, the data collected in a parameter matrix $\boldsymbol{N}_{k}^{z}$ provides an alternative approach
for building local ROBs by a mechanism centered around updating and splitting. Again, the updating criterion,
which is proposed in this work for avoiding extrapolation, is described in Algorithm \ref{alg:criterion}.
Its basic principle is described in Algorithm \ref{alg:update}. The proposed splitting algorithm is described
below.

\vglue 0.125truein
\begin{algorithm}[H] 
  \DontPrintSemicolon
  \KwData{Sampled parameter vector $\boldsymbol{z}_{k}^{\star}$ (or point in $\mathcal{D}_k$), corresponding
          HDM-based solution snapshot $\boldsymbol{u}_{k}^{\star}$, ROB $\boldsymbol{V}_{k}$, matrix
          of sampled deformation gradient parameter points $\boldsymbol{N}_{k}^{z}$, and snapshot matrix
          $\boldsymbol{N}_{k}^{u}$} 
  \KwResult{Updated ROB $\boldsymbol{V}_{k}$, updated sampled parameter matrix $\boldsymbol{N}_{k}^{z}$,
            and updated snapshot matrix $\boldsymbol{N}_{k}^{u}$} 
  Update $\boldsymbol{N}_{k}^{z} \leftarrow \left [\boldsymbol{N}_{k}^{z} \quad \boldsymbol{z}_{k}^{\star}\right ]$\; 
  Update $\boldsymbol{N}_{k}^{u} \leftarrow \left [\boldsymbol{N}_{k}^{u} \quad \boldsymbol{u}_{k}^{\star}\right ]$\;
  Perform low-rank SVD update of $\boldsymbol{N}_{k}^{u}$\;
  Extract new (updated) ROB $\boldsymbol{V}_k$\;
  \caption{Updating a ROB.} 
  \label{alg:update} 
\end{algorithm}
\vglue 0.125truein

Specifically, $\boldsymbol{N}_{k}^{z}$ provides a mean for splitting a region of the parameter space
$\mathcal{D}_k$ rather than a region of the solution space in two parts, using for example Principle
Component Analysis (PCA) or POD, both of which build on SVD. PCA is most useful for exploratory data
analysis. It is used here to compute the direction of maximum variance $\boldsymbol{n}_{k}$ of a set
of sampled parameter points ${\cal S}_z^{m_k}$ as the first eigenvector of the covariance of the corresponding
parameter matrix
\begin{equation}
\label{eq:cov}
  \text{Cov}(\boldsymbol{N}_{k}^{z})=\boldsymbol{N}_{k}^{z}\boldsymbol{N}_{k}^{z^T}-\bar{\boldsymbol{z}}_{k}\bar{\boldsymbol{z}}_{k}^{T}
\end{equation}
where $\bar{\boldsymbol{z}}_{k} \in \mathbb{R}^{n_{\mathrm{sd}}}$ is the columnwise mean of $\boldsymbol{N}_{k}^{z}$
-- that is,
\begin{equation}
  \bar{\boldsymbol{z}}_{k_i} = \displaystyle{\frac{1}{m_k}\sum_{j=1}^{m_k}\boldsymbol{N}_{k_{ij}}^{z}, \qquad i=1,\ldots,n_{\mathrm{sd}}}
\label{eq:CWCG}
\end{equation}
and therefore $\bar{\boldsymbol{z}}_{k}$ represents the algebraic centroid of the sampled parameter points.

\begin{figure}[h!]
  \centering
  \begin{subfigure}[b]{0.45\textwidth}
    \centering
    \includegraphics[width=\textwidth]{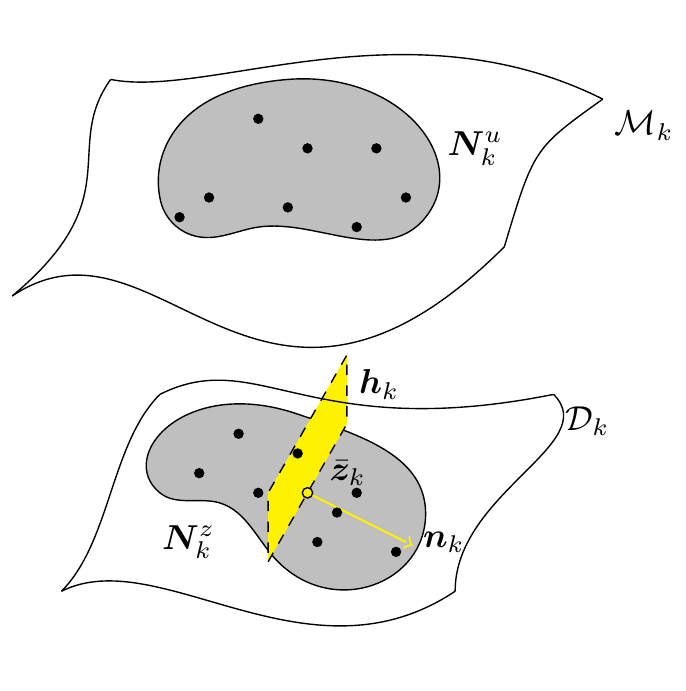}
    \captionsetup{justification=centering}
    \caption{splitting a set of sampled parameter points into two sets using a hyperplane}
    \label{fig:split1}
  \end{subfigure}
  \begin{subfigure}[b]{0.45\textwidth}
    \centering
    \includegraphics[width=\textwidth]{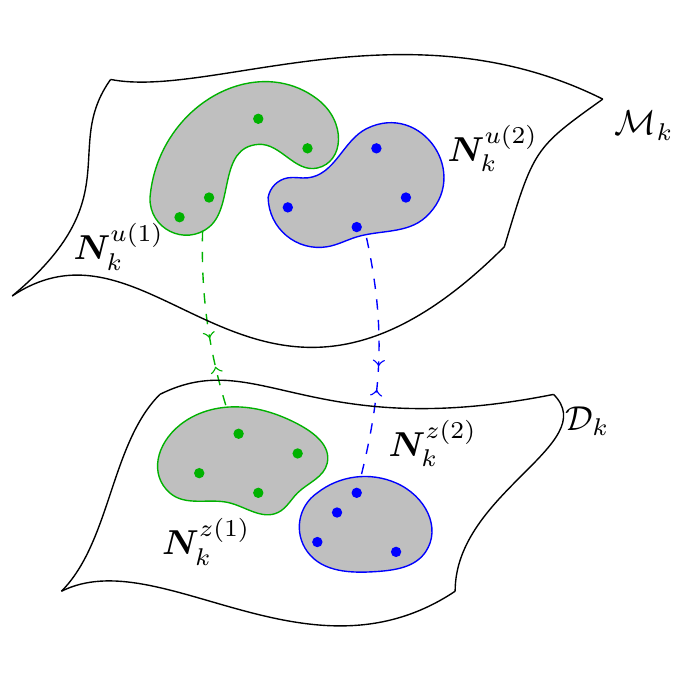}
    \captionsetup{justification=centering}
    \caption{partitioning the corresponding snapshot matrix into two matrices accordingly}
    \label{fig:split2}
  \end{subfigure}
  \caption{Splitting a set of sampled parameter points in two sets in order to partition a snapshot matrix
           into two matrices and thereby split a ROB into two local ROBs.}
  \label{fig:split}
\end{figure}

Next, the hyperplane $\boldsymbol{h}_{k}$ containing the centroid point $\bar{\boldsymbol{z}}_{k}$ and
orthogonal to $\boldsymbol{n}_{k}$ is constructed. This hyperplane can be described as
\begin{equation}
\label{eq:plane}
  \boldsymbol{h}_{k} = \left\{\boldsymbol{x} | (\boldsymbol{x}-\bar{\boldsymbol{z}}_{k})^{T}\boldsymbol{n}_{k}=0\right\}
\end{equation}
It naturally splits ${\cal S}_z^{m_k}$ into two smaller sets, where each partitioned set has a more compact
spatial distribution in $\mathcal{D}_k$ (see Figure \ref{fig:split2}). 

As illustrated in Figure \ref{fig:split}, the splitting of a set of sampled parameter points ${\cal S}_z^{m_k}$
into two smaller sets, which is equivalent to the partitioning of a matrix $\boldsymbol{N}_{k}^{z}$ into
two smaller matrices, implies a splitting of a region of the parameter space $\mathcal{D}_k$. It also
leads to the partitioning of the snapshot matrix $\boldsymbol{N}_{k}^{u}$ in two snapshot matrices. Finally,
the compression of these two matrices leads in turn to the construction of two new local ROBs. The entire
splitting process is summarized in Algorithm \ref{alg:split}.

\begin{algorithm}[H]
\DontPrintSemicolon
  \KwData{Matrix of sampled deformation gradient parameter points $\boldsymbol{N}_{k}^{z}$, snapshot
          matrix $\boldsymbol{N}_{k}^{u}$, and database of HPROM-related information $\mathcal{DB}_k$} 
  \KwResult{Split local ROBs $\boldsymbol{V}_{k}^{(1)}$ and $\boldsymbol{V}_{k}^{(2)}$, split parameter
            matrices $\boldsymbol{N}_{k}^{z^{(1)}}$ and $\boldsymbol{N}_{k}^{z^{(2)}}$, split snapshot 
            matrices $\boldsymbol{N}_{k}^{u^{(1)}}$, and $\boldsymbol{N}_{k}^{u^{(2)}}$, and updated
            database $\mathcal{DB}_k$}
  Compute the first eigenvector $\boldsymbol{n}_{k}$ of the matrix $\text{Cov}(\boldsymbol{N}_{k}^{z})$ \eqref{eq:cov}\;
  Construct the hyperplane $\boldsymbol{h}_{k}$ \eqref{eq:plane}\;
  Using $\boldsymbol{h}_{k}$, partition $\boldsymbol{N}_{k}^{z}$ into two matrices $\boldsymbol{N}_{k}^{z^{(1)}}$
  and $\boldsymbol{N}_{k}^{z^{(2)}}$\; 
  Partition $\boldsymbol{N}_{k}^{u}$ accordingly, into two matrices $\boldsymbol{N}_{k}^{u^{(1)}}$ and
  $\boldsymbol{N}_{k}^{u^{(2)}}$\;
  Compute $\boldsymbol{V}_{k}^{(1)}$, $\boldsymbol{V}_{k}^{(2)}$ from $\boldsymbol{N}_{k}^{u^{(1)}}$,
  $\boldsymbol{N}_{k}^{u^{(2)}}$, respectively\; 
  Add $\boldsymbol{V}_{k}^{(1)}$, $\boldsymbol{N}_{k}^{z^{(1)}}$, $\boldsymbol{N}_{k}^{u^{(1)}}$, and
  $\boldsymbol{V}_{k}^{(2)}$, $\boldsymbol{N}_{k}^{z^{(2)}}$, $\boldsymbol{N}_{k}^{u^{(2)}}$ to the content
  of the database $\mathcal{DB}_k$\;
  Delete $\boldsymbol{V}_{k}$, $\boldsymbol{N}_{k}^{z}$, $\boldsymbol{N}_{k}^{u}$ from the content of
  the database $\mathcal{DB}_k$ 
  \caption{Splitting a ROB.} 
  \label{alg:split}
\end{algorithm}

\subsection{Approximation using a most-appropriate local ROB}

The adaptive, in-situ sampling process described so far, which centers around updating and splitting
ROBs, eventually leads to the presence of multiple local ROBs in a database $\mathcal{DB}_k$. Hence,
when a microscale BVP characterized by an unsampled parameter point $\boldsymbol{z}_{k}^{\star}$ is queried
for fast solution by a local HPROM, the most-appropriate local ROB available in $\mathcal{DB}_k$ must be
identified online. This can be performed as proposed below.

Let $\mathcal{S}_V^{n_k} = \left\{\boldsymbol{V}_{k}^{(1)}, \ldots,\boldsymbol{V}_{k}^{(n_k)}\right\}$
denote the set of all $n_k$ local ROBs that are pre-computed and stored at a given time in a database
$\mathcal{DB}_k$. Recall from Section \ref{sec:adaptive-ROM} that each HDM-based solution $\boldsymbol{u}_{k}^{(i)}$
of a microscale BVP is associated with a parameter point $\boldsymbol{z}_{k}^{(i)}$ of the parameter
space $\mathcal{D}_k$ that characterizes the BVP. Hence, a split snapshot matrix $\boldsymbol{N}_{k}^{u^{(l)}}$
can be associated on the same basis with a split parameter matrix $\boldsymbol{N}_{k}^{z^{(l)}}$. Similarly,
a local ROB $\boldsymbol{V}_{k}^{(l)}$ can be associated with the snapshot matrix $\boldsymbol{N}_{k}^{u^{(l)}}$
that was compressed to obtain $\boldsymbol{V}_{k}^{(l)}$, and by transitivity, with the split parameter
matrix $\boldsymbol{N}_{k}^{z^{(l)}}$.

Let $\bar{\boldsymbol{z}}_{k}^{(l)}$ denote the algebraic centroid of the parameter matrix $\boldsymbol{N}_{k}^{z^{(l)}}$
(defined as in \eqref{eq:CWCG}). Following the same reasoning as above, it follows that $\bar{\boldsymbol{z}}_{k}^{(l)}$
can be associated with the local ROB $\boldsymbol{V}_{k}^{(l)}$. More generally, a bijection can be defined
between the set of sampled parameter points $\bar{\boldsymbol{z}}_{k}^{(l)}$ and the set of pre-computed
local ROBs $\boldsymbol{V}_{k}^{(l)}$. 

Now, given a queried BVP characterized by an unsampled deformation gradient parameter point $\boldsymbol{z}_{k}^{\star}$
-- or more simply, given a queried but unsampled parameter point $\boldsymbol{z}_{k}^{\star}$ -- a most-appropriate,
pre-computed, local ROB $\boldsymbol{V}_{k}^{(j_k)} \in \mathcal{S}_V^{n_k}$ -- that is, $j_k \in \left\{1, \ldots, n_k\right\}$,
can be defined as the pre-computed local ROB whose corresponding parameter point is the closest to the
queried but unsampled parameter point $\bar{\boldsymbol{z}}_{k}^{\star}$ -- that is,
\begin{equation*}
\label{eq:search-ROB}
  j_{k}=\underset{i \in \left\{1,\ldots,n_k\right\}}{\operatorname{argmin}} \quad \left\|\bar{\boldsymbol{z}}_{k}^{(i)}
  -\boldsymbol{z}_{k}^{\star}\right\|
\end{equation*}
where $\|~\|$ denotes any preferred norm, for example, $\|~\|_2$.

\subsection{Summary: in-situ adaptive microscale PMOR framework}

Algorithm \ref{alg:summary} summarizes the in-situ, adaptive, PMOR framework proposed for the acceleration
of the solution of the microscale BVPs at any level $k$, $k=1,\ldots,n_s$.

\begin{algorithm}[H]
\DontPrintSemicolon 
  \KwData{User-specified residual tolerance $r_{\text{tol}}$, and user-specified ROB capacity limit $c_{\text{max}}$} 
  \KwResult{HDM/HPROM-based solutions of $\Pi_{j=0}^{k-1}\left|\mathcal{G}_{j}\right|$ microscale BVPs
            at level $k$, and database of HPROM-related information $\mathcal{DB}_k$}
  Initialize $m_k$\;
  Collect $\boldsymbol{N}_{k}^{z}$ and $\boldsymbol{N}_{k}^{u}$ and compute $\boldsymbol{V}_{k}$ using
  the in-situ training strategy described in Section \ref{sec:rom-micro}\; 
  Initialize the content of $\mathcal{DB}_k$ with $\boldsymbol{V}_{k}$, $\boldsymbol{N}_{k}^{z}$ and
  $\boldsymbol{N}_{k}^{u}$\; 
  \For{solving each of the remaining $\prod_{j=0}^{k-1}\left|\mathcal{G}_{j}\right|-m_{k}$ microscale
  BVPs}{
    Compute the deformation gradient parameter point $\boldsymbol{z}_{k}^{\star}$\; 
    Identify in $\mathcal{DB}_k$ the corresponding most-appropriate local ROB $\boldsymbol{V}_{k}^{(j_k)}$\; 
    Apply Algorithm \ref{alg:criterion} to update or not $\mathcal{DB}_k$\; 
    \eIf{\textrm{UF} = \textrm{true} (see Algorithm \ref{alg:criterion})}{
      Update $\mathcal{DB}_k$ using Algorithm \ref{alg:update}\; 
      \eIf{size of updated ROB $> c_{\text{max}}$}{ 
        Split ROB using Algorithm \ref{alg:split}
      }{ 
        Continue (to the next microscale BVP)\;
      }
    }{ 
      Continue (to the next microscale BVP)\; 
    }
  }
  \caption{In-situ, adaptive, PMOR framework for the microscale BVPs at level $k$.} 
  \label{alg:summary} 
\end{algorithm}

\section{Applications}
\label{sec:apps}

Here, the computational feasibility as well as the performance of the in-situ, adaptive, PMOR framework
for microscale BVPs are demonstrated for two nonlinear, nonparametric, dynamic, two-scale simulations
($n_s = 1$).  Specifically, this framework is applied to the solution of two academic but nevertheless
computationally intensive, nonlinear, two-level, dynamic response problems. For both applications, all
ROBs are constructed using POD based on displacement solution snapshots. Because adaptive hyperreduction
is beyond the scope of this paper, all HPROMs are constructed based on a single reduced mesh that is
trained and generated using ECSW during the initialization step of Algorithm \ref{alg:summary}. Given
that at the macroscale level neither considered problem is parametric, PMOR is applied in both cases
only to the microscale model. This setting has two benefits:
\begin{itemize}
  \item It enables the exclusive focusing on the performance of the proposed in-situ, adaptive approach
        for performing the PMOR of microscale problems.
  \item It corresponds to the worst case scenario of PMOR performance, as any computational overhead
        incurred by the PMOR process may only be amortized in this case by the same simulation to which
        the PMOR process is applied. Hence, any reported speedup factor is a genuine speedup factor.
\end{itemize}
		
Since both applications considered herein involve two-scale simulations ($n_s = 1$), the model labels
PROM and HPROM refer throughout the remainder of this section to the microscale computational model.
On the other hand, the label HDM may refer to either the macroscale or the microscale computational model.

In both applications, the tolerance for the residual-based error indicator is set to $r_{\text{tol}}=0.001$,
and the capacity limit for the dimension of a local ROB is set to $c_{\text{max}}=20$. The accuracy of
all PROM and HPROM approximations is evaluated using the following definition of a relative global error
\begin{equation*}
  e_{\text{rel}}=\left|\frac{\sum\limits_{t \in P} \Delta s\left(\boldsymbol{u}_{\diamond}(t) -\tilde{\boldsymbol{u}}_{\diamond}(t)\right)^{T}
  \boldsymbol{1}}{\sum\limits_{t \in P} \Delta s \boldsymbol{u}_{\diamond}^{T}(t) \boldsymbol{1}}\right| \times 100 \%
\end{equation*}
where the $\diamond$ subscript designates the $x$, $y$, or $z$ direction of a global reference frame,
$\boldsymbol{u}_{\diamond}(t) \in \mathbb{R}^{n_0^d}$ is the vector of $\diamond$-displacements at time
$t$ extracted from the macroscale solution based on HDMs at both the macroscale and microscale levels,
$\tilde{\boldsymbol{u}}_{\diamond}(t) \in \mathbb{R}^{n_0^d}$ is the vector of $\diamond$-displacements
at time $t$ extracted from the macroscale solution but based on a PROM or HPROM at the microscale level,
$\Delta s$ is the error sampling time-step size, $\boldsymbol{1}$ denotes the vector of 1 entries, and
$P$ is the set of time-stamps
\begin{equation*}
  P=\left\{t_{0}, t_{0}+\Delta s, t_{0}+2 \Delta s, \ldots, t_{f}\right\}
\end{equation*}
where $t_{0}$ is the initial simulated time and $t_{f}$ is the final simulated time.

All computations reported herein are performed on a Linux cluster with 84 compute nodes, dual-socket
Intel Xeon E5-2650 v2 CPUs (8 cores/socket), 64 GB of DDR3 memory/node, a 1:1 subscribed FDR Infiniband
network, and a Lustre shared file system tuned for parallel I/O. 

For both applications, all reported PROM and HPROM wall-clock timings as well as speedup factors account
for all training costs.

\subsection{Vibration of a nonlinear solid elastic bar with a microvoid structure}

A computational model for a solid bar made of a heterogeneous structure is considered here. At the finest
scale, this structure is characterized by microscopic voids and a porosity ratio of $0.5$ (here, the
porosity ratio is defined as the ratio of the total volume occupied by the voids, and the total volume
occupied by the material and the voids). The bar has a solid square cross section of height and width
$H=1$ cm, and a length $L=20$ cm. It is fixed at both ends. A distributed external force generated by
a uniform pressure of magnitude $p=100$ MPa is suddenly applied on its upper face and oriented initially
in the transverse direction. Due to symmetry, only half of the bar along its length is modeled.

\begin{figure}[h]
\centering
\includegraphics[width=0.8\textwidth]{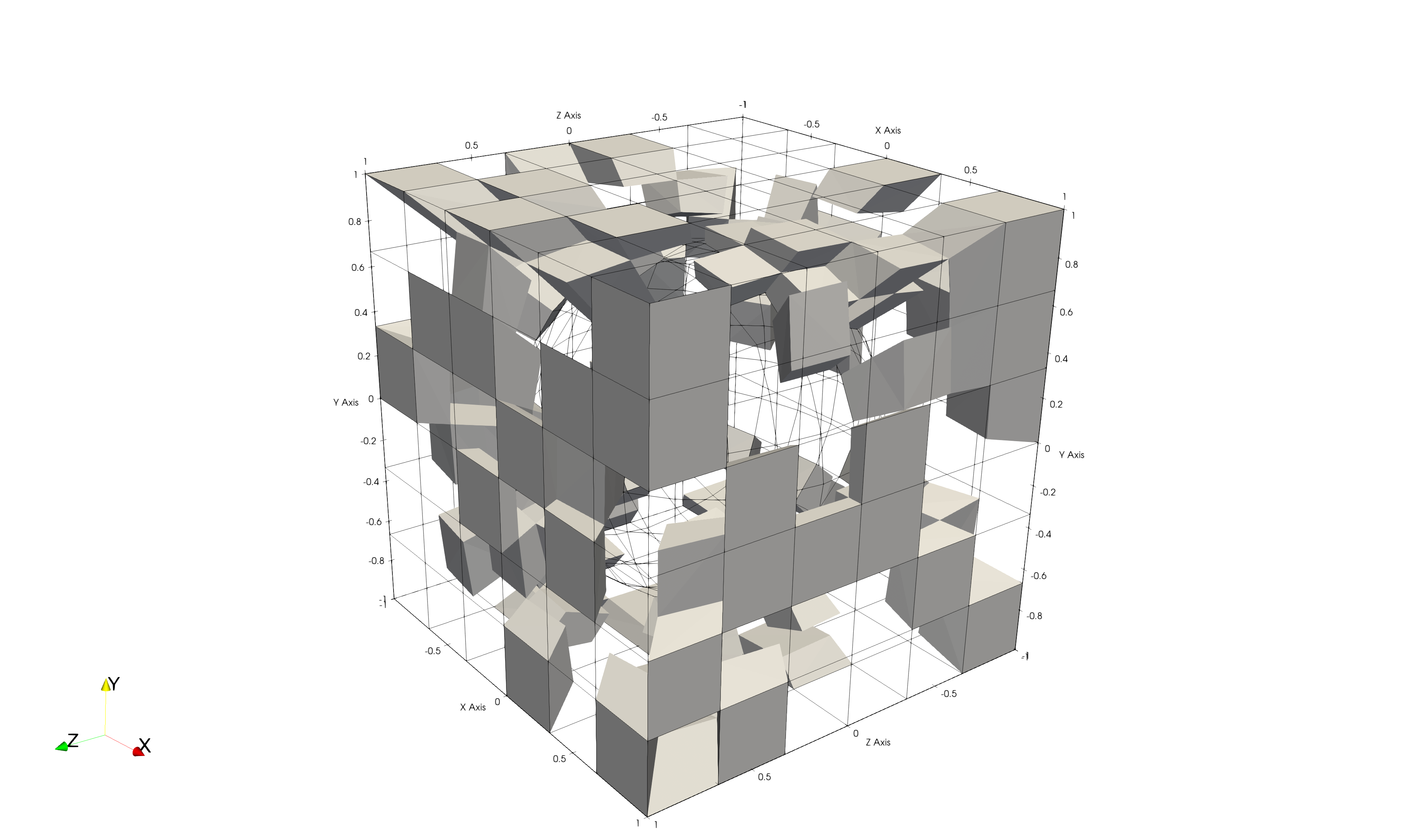}
\caption{Vibration of a nonlinear solid elastic bar with a microvoid structure: FE mesh underlying the
         microscale HDM (wireframes) and associated ECSW-generated reduced mesh (shaded elements).}
\label{fig:p1_micro_mesh}
\end{figure}

At the microscale level, a FE model of the bar (see Figure \ref{fig:p1_micro_mesh}) is constructed using
648 eight-noded hexahedral elements with three displacement DOFs per node. The dimension of this HDM
is $n_1^d = 1,962$ (1,962 DOFs). At this level, the material response is assumed to be well represented
by a compressible neo-Hookean (finite strains) hyperelastic constitutive law with Young's modulus $E=207$
GPa and Poisson's ratio $\nu=0.3$. Uniform essential boundary conditions are adopted at this level.

\begin{figure}[h]
\centering
\includegraphics[width=0.8\textwidth]{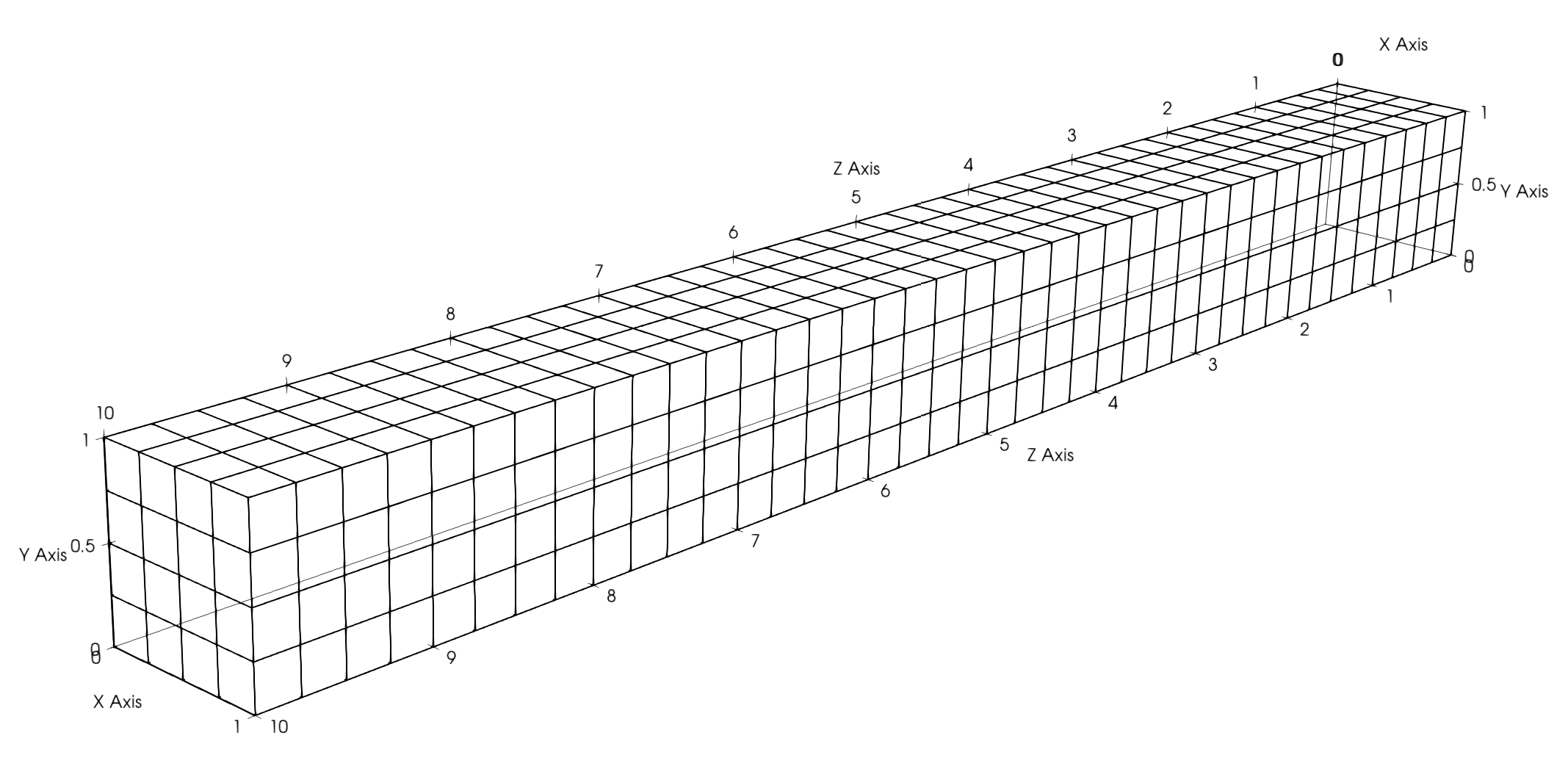}
\caption{Vibration of a nonlinear solid elastic bar with a microvoid structure: FE mesh underlying the
         macroscale HDM.}
\label{fig:p1_macro_mesh}
\end{figure}

At the macroscale level, a FE model of the bar (see Figure \ref{fig:p1_macro_mesh}) is constructed using
640 eight-noded hexahedral elements with three displacement DOFs/node, and a total number of $2,975$ DOFs
($n_0^d = 2,975$). At this level, the material density is $\rho_{0}=7,830$ kg/m$^3$, and temporal discretization
is performed using the explicit central difference scheme and the constant time-step $\Delta t=3.25 \times 10^{-7}$ s.
The response of the bar is computed from $t_{0}=0$ s until $t_{f}=1 \times 10^{-2}$ s.

First, an initial nonlinear PROM is constructed for the microscale computational model using the solution
snapshots collected at $t_{0}=0$ s during the HDM-based solution of the first $m_1=500$ microscale BVPs
encountered at the Gauss points of the macroscale computational model. Then, this PROM is hyperreduced
using ECSW configured with the Lawson and Hanson sparse NNLS algorithm and a relative accuracy threshold
$\varepsilon=10^{-3}$, which generates a reduced mesh with 115 elements (${\widetilde{\cal E}} = 115$).

Several simulations are performed for this problem using: the macroscale and microscale HDMs described
above; the proposed in-situ, adaptive, microscale PMOR framework, with and without hyperreduction (the
reader is reminded that in this work, any reduced mesh generated by ECSW-based hyperreduction is not
updated); and the counterpart nonadaptive microscale PMOR framework described in \cite{zahr2017multilevel},
with and without hyperreduction.

Table \ref{tab:adaptive} reports the speedup factors and relative global errors computed using the error
sampling time-step $\Delta s=\Delta t$ measured for the proposed in-situ, adaptive, microscale PMOR with
and without hyperreduction. Similarly, Tables \ref{tab:non1}--\ref{tab:non3} report these performance
results for the counterpart in-situ, nonadaptive, microscale PMOR framework previously proposed in \cite{zahr2017multilevel}.
Since the latter framework depends in this case on the ad-hoc parameter $n_{0}^{\mathcal G}$ specifying
the number of microscale BVPs to solve in order to train a global microscale PROM and the corresponding
HPROM, three different values of this parameter are considered. The first one, $n_{0}^{\mathcal{G}} = 62$,
corresponds to 0.2\% of the total number $n_{ms}$ of microscale problems solved during the multiscale
simulation performed in the time-interval $[0, 1 \times 10^{-2}]$ s using the time-step size $\Delta t = 3.25 \times 10^{-7}$ s.
The second and third values, $n_{0}^{\mathcal{G}} = 154$ and $n_{0}^{\mathcal{G}} = 308$, correspond
to 0.5\% of $n_{ms}$ and 1\% of $n_{ms}$, respectively. From the results reported in this table, the
following observations are noteworthy:
\begin{itemize}
  \item In all cases, the HPROM delivers a much better speedup factor than its underlying PROM. This
        is expected given that the problem considered here is geometrically nonlinear (finite strains).
  \item The speedup factors delivered by all HPROMs are impressive, considering that the macroscale model
        is not even reduced.
  \item The adaptively computed PROM delivers 6 to 50 times more accurate results than its nonadaptive
         counterparts trained with any considered value of $n_{0}^{\mathcal{G}}$, and performs about
         1.25 times faster.
  \item The adaptively computed HPROM -- except for the nonadapted but fixed reduced mesh -- delivers
        2 to 3 times more accurate results than its nonadaptive counterpart trained with $n_{0}^{\mathcal{G}} = 62$,
        but is 1.35 times slower. It delivers roughly the same wall-clock time performance as its nonadaptive
        counterpart trained with $n_{0}^{\mathcal{G}} = 154$, and about the same accuracy. On the other
        hand, it is about 1.2 times faster than its nonadaptive counterpart trained with $n_{0}^{\mathcal{G}} = 308$,
        while delivering essentially the same accuracy.
\end{itemize}
It follows that in all considered cases, the adaptive PROM is faster than its nonadaptive counterpart
and considerably more accurate. On the other hand, because the reduced mesh resulting from the ECSW-based
hyperreduction is fixed in this work, the otherwise adaptive HPROM can be anywhere from slightly slower
to slightly faster than its nonadaptive counterpart, and roughly as accurate.

\begin{figure}[h!]
\centering
\begin{subfigure}[b]{0.48\textwidth}
\centering
\includegraphics[width=\textwidth]{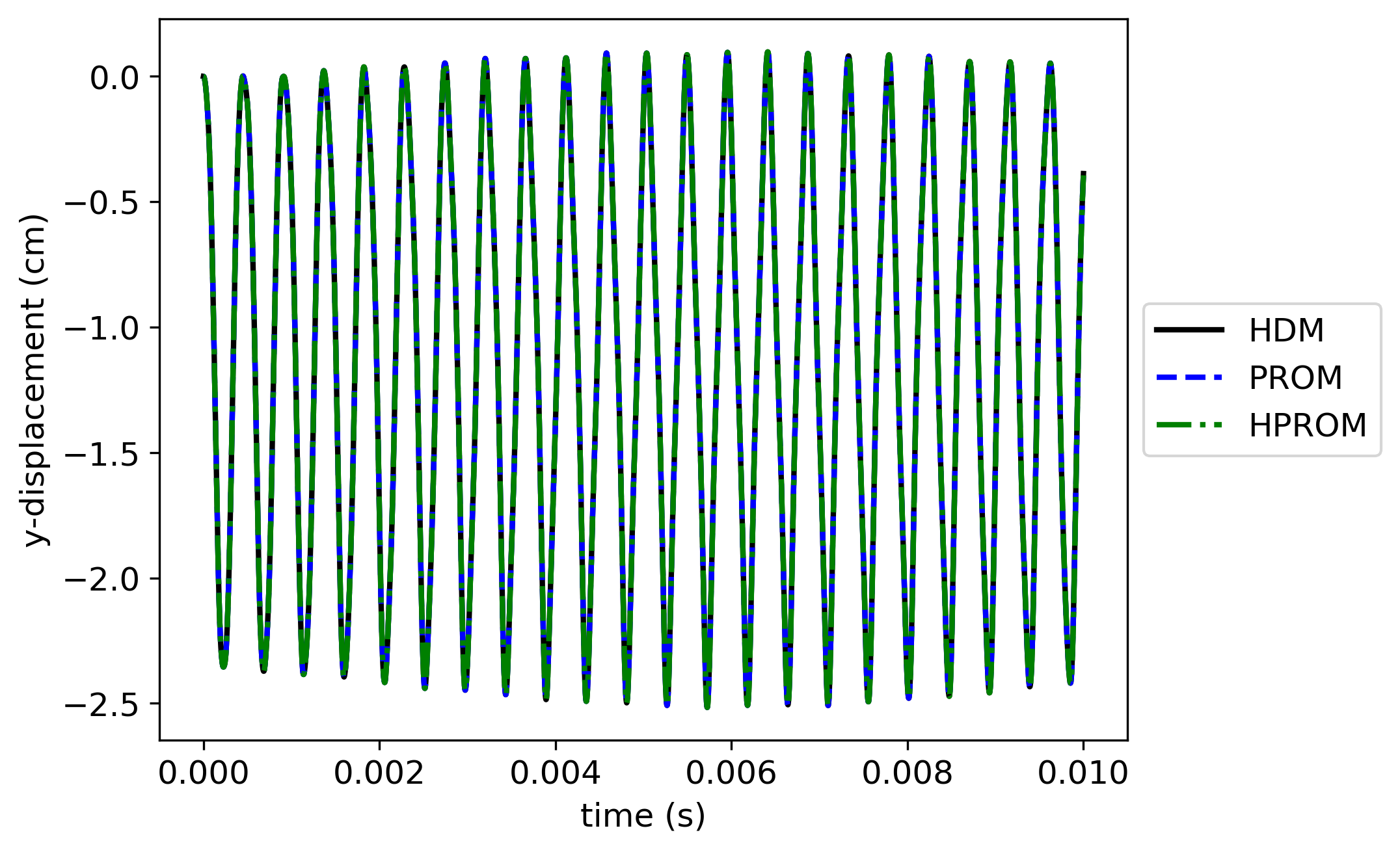}
\captionsetup{justification=centering}
\caption{displacement}
\label{fig:p1_y_disp}
\end{subfigure}
\begin{subfigure}[b]{0.48\textwidth}
\centering
\includegraphics[width=\textwidth]{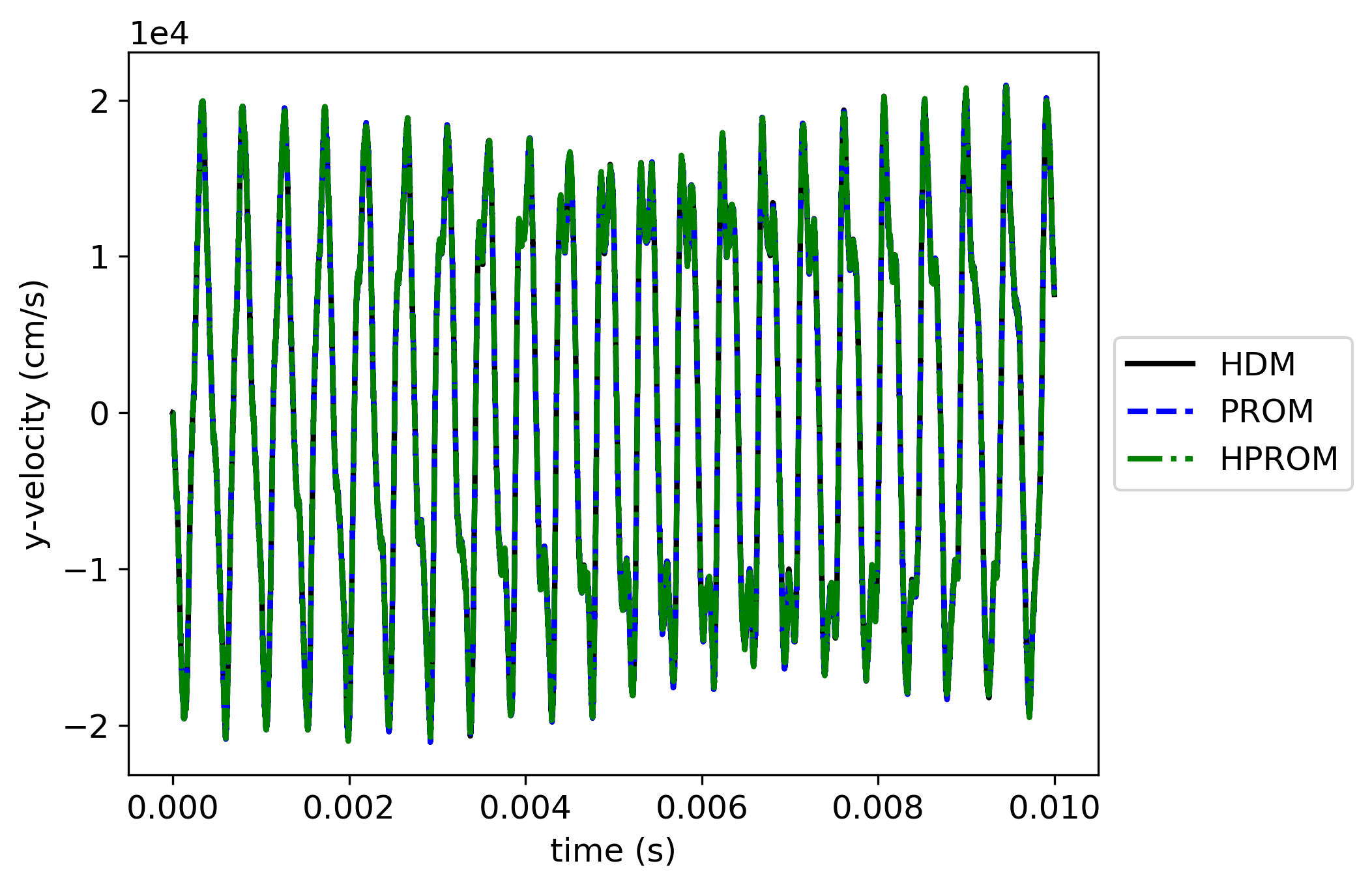}
\captionsetup{justification=centering}
\caption{velocity}
\label{fig:p1_y_velo}
\end{subfigure}
\caption{Vibration of a nonlinear solid elastic bar with a microvoid structure -- time-histories of a
         displacement DOF and a velocity DOF at a node of the mesh underlying the macroscale HDM computed 
         using: the macroscale and microscale HDMs; and the in-situ, adaptive, microscale framework for
         PMOR with and without hyperreduction.}
\label{fig:p1_history}
\end{figure}

Considering that all aforementioned performance results account for all computational overhead, including
the adaptation overhead, and that the proposed in-situ, adaptive, microscale PMOR framework is fully
automated, these results demonstrate not only the feasibility of this framework but also its superiority
over its nonadaptive counterpart. In particular, Figure \ref{fig:p1_history} shows that even when the
reduced mesh underlying the microscale computational model is not adapted, the proposed adaptive framework
delivers a remarkable accuracy.

\begin{table}[h!]
\centering
\captionsetup{justification=centering}
\caption{Vibration of a nonlinear solid elastic bar with a microvoid structure: performance assessment.}
\begin{subtable}[h!]{0.45\textwidth}
\begin{tabular}{ccc}
\hiderowcolors
\hline 
Model & \shortstack{Wall-clock \\ time (s)} & \shortstack{Speedup \\ factor}\\
\hline 
HDM & $1.35 \times 10^5$ & N.A.\\
PROM & $3.75 \times 10^4$ & $3.60$ \\
HPROM & $8.92 \times 10^3$ & $15.2$ \\
\hline 
Error measure $e$ & PROM & HPROM\\
\hline 
y-displacement (\%) & $0.0227$ & $0.440$\\
z-displacement (\%) & $0.0435$ & $0.747$\\
y-velocity (\%) & $0.103$ & $1.77$\\
z-velocity (\%) & $0.986$ & $6.68$\\
\hline
\end{tabular}
\centering
\captionsetup{justification=centering}
\caption{Adaptive}
\label{tab:adaptive}
\end{subtable}
\hfill
\begin{subtable}[h!]{0.45\textwidth}
\centering
\begin{tabular}{ccc}
\hiderowcolors
\hline 
Model & \shortstack{Wall-clock \\ time (s)} & \shortstack{Speedup \\ factor}\\
\hline 
HDM & $1.35 \times 10^5$ & N.A.\\
PROM & $4.64 \times 10^4$ & $3.03$ \\
HPROM & $6.58 \times 10^3$ & $20.5$ \\
\hline 
Error measure $e$ & PROM & HPROM\\
\hline 
y-displacement (\%) & $1.23$ & $1.27$\\
z-displacement (\%) & $2.11$ & $2.59$\\
y-velocity (\%) & $3.05$ & $3.23$\\
z-velocity (\%) & $9.52$ & $10.6$\\
\hline
\end{tabular}
\captionsetup{justification=centering}
	\caption{Nonadaptive ($n_{0}^{\mathcal{G}} = 62$)}
\label{tab:non1}
\end{subtable}
\hfill
\begin{subtable}[h!]{0.45\textwidth}
\centering
\begin{tabular}{ccc}
\hiderowcolors
\hline 
Model & \shortstack{Wall-clock \\ time (s)} & \shortstack{Speedup \\ factor}\\
\hline 
HDM & $1.35 \times 10^5$ & N.A.\\
PROM & $4.54 \times 10^4$ & $2.97$ \\
HPROM & $8.71 \times 10^3$ & $15.5$ \\
\hline 
Error measure $e$ & PROM & HPROM\\
\hline 
y-displacement (\%) & $0.0530$ & $0.134$\\
z-displacement (\%) & $0.106$ & $0.219$\\
y-velocity (\%) & $1.11$ & $1.21$\\
z-velocity (\%) & $6.24$ & $6.50$\\
\hline
\end{tabular}
\captionsetup{justification=centering}
	\caption{Nonadaptive ($n_{0}^{\mathcal{G}} = 154$)}
\label{tab:non2}
\end{subtable}
\hfill
\begin{subtable}[h]{0.45\textwidth}
\centering
\begin{tabular}{ccc}
\hiderowcolors
\hline 
Model & \shortstack{Wall-clock \\ time (s)} & \shortstack{Speedup \\ factor}\\
\hline 
HDM & $1.35 \times 10^5$ & N.A.\\
PROM & $4.68 \times 10^4$ & $2.88$ \\
HPROM & $1.06 \times 10^4$ & $12.7$ \\
\hline 
Error measure $e$ & PROM & HPROM\\
\hline 
y-displacement (\%) & $0.0271$ & $0.108$\\
z-displacement (\%) & $0.0762$ & $0.155$\\
y-velocity (\%) & $1.10$ & $1.03$\\
z-velocity (\%) & $6.07$ & $5.79$\\
\hline
\end{tabular}
\captionsetup{justification=centering}
	\caption{Nonadaptive ($n_{0}^{\mathcal{G}} = 308$)}
\label{tab:non3}
\end{subtable}
\label{tab:p1}
\end{table}

\newpage
\subsection{Impact of a nonlinear elastic cylinder with a composite microstructure on a rigid wall}

Next, a computational model is considered for a cylindrical projectile with a heterogeneous microscale
structure characterized by a stiff, unidirectional fiber embedded in a flexible matrix. The cylinder
has a radius $R=0.476$ cm and a length $L=1.27$ cm. It has an initial velocity $v=15,000$ cm$/$s in the
negative $z$-direction, and a frictionless rigid wall boundary condition at $z=0$. Due to symmetry, only
one quarter of the bar along its length is modeled. This problem is a variant of the so-called Taylor
impact problem \cite{taylor1948use}.

For this problem, a FE microscale computational model is constructed using $729$ eight-noded hexahedral
elements with three displacement DOFs per node, and a total number of $2,988$ DOFs. Uniform essential
boundary conditions are adopted at this level. The microscale material is modeled as follows:
\begin{itemize} 
  \item The fiber is represented as a compressible neo-Hookean (finite strains) material with Young's
        modulus $E=212.52 \times 10^{10}$ dyne$/$cm$^2$ and Poisson's ratio $\nu=0.4$. It runs parallel
        to the local $x$-axis along the center of the cube where the locally attached microstructure
        is defined. Its cross section is approximated by a square of edge size equal to one-third of
        the unit cell width.  
  \item The flexible matrix surrounding this fiber is also represented as a compressible neo-Hookean
        material (finite strains), but with Young's modulus  $E=72.52 \times 10^{10}$ dyne$/$cm$^2$ and
        Poisson's ratio $\nu=0.4$.
\end{itemize}

\begin{figure}[h!]
\centering
\includegraphics[width=1.0\textwidth]{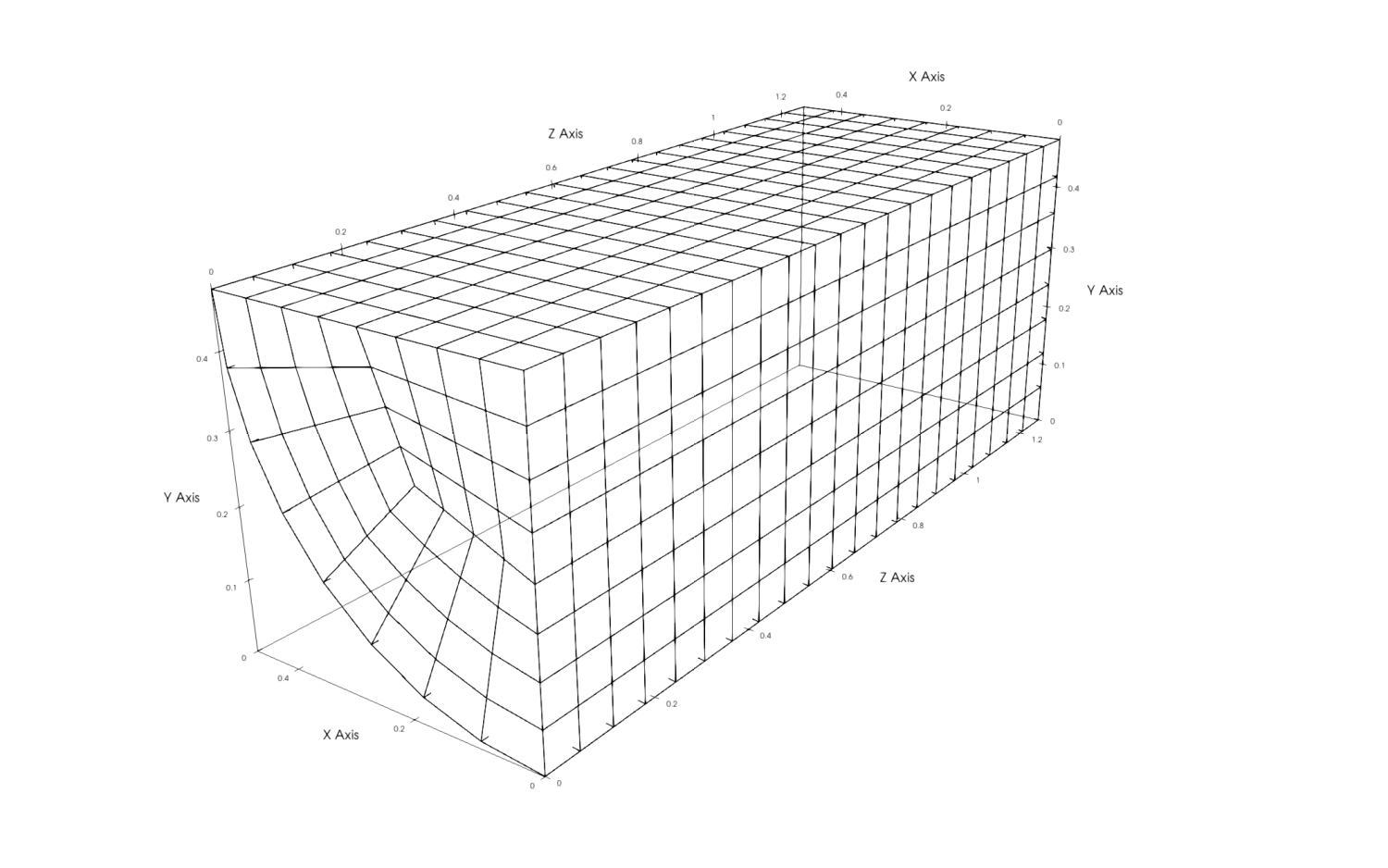}
\caption{Impact on a rigid wall of a nonlinear elastic cylinder with a composite microstructure: FE mesh
         underlying the macroscale HDM.}
\label{fig:p2_macro_mesh}
\end{figure}

At the macroscale level, a FE model is constructed using $1,008$ eight-noded hexahedral elements with
three displacement DOFs/node, and a total number of $1,342$ DOFs (see Figure \ref{fig:p2_macro_mesh}).
At this level, the density is $\rho_{0}=2.7$ g/cm$^3$. Temporal discretization is performed using the
explicit central difference scheme and the constant time-step $\Delta t= 1.0 \times 10^{-9}$ s. For this
problem, the multiscale simulation is performed from $t_{0}=0$ s until $t_{f}=1 \times 10^{-6}$ s.

As in the previous application, the proposed in-situ, adaptive, microscale framework for PMOR is initialized
with the construction of a nonlinear PROM based on the first $m_1=500$ HDM-based solutions of the microscale
BVPs encountered at the Gauss points of the macroscale computational model. Then, this PROM is hyperreduced
using ECSW configured with the same sparse NNLS algorithm as in the previous application, and the same
threshold for relative accuracy $\varepsilon=10^{-3}$. ECSW generates in this case a reduced mesh with
88 elements (${\widetilde{\cal E}} = 88$). Three multiscale computations are performed using: the aforementioned
macroscale and microscale HDMs; and the proposed in-situ, adaptive, microscale PMOR framework, with and
without hyperreduction. The obtained performance results are reported in Table \ref{tab:p2}, where the
relative global errors are computed with $\Delta s=\Delta t$.

\begin{figure}[h!]
\centering
\begin{subfigure}[b]{0.48\textwidth}
\centering
\includegraphics[width=\textwidth]{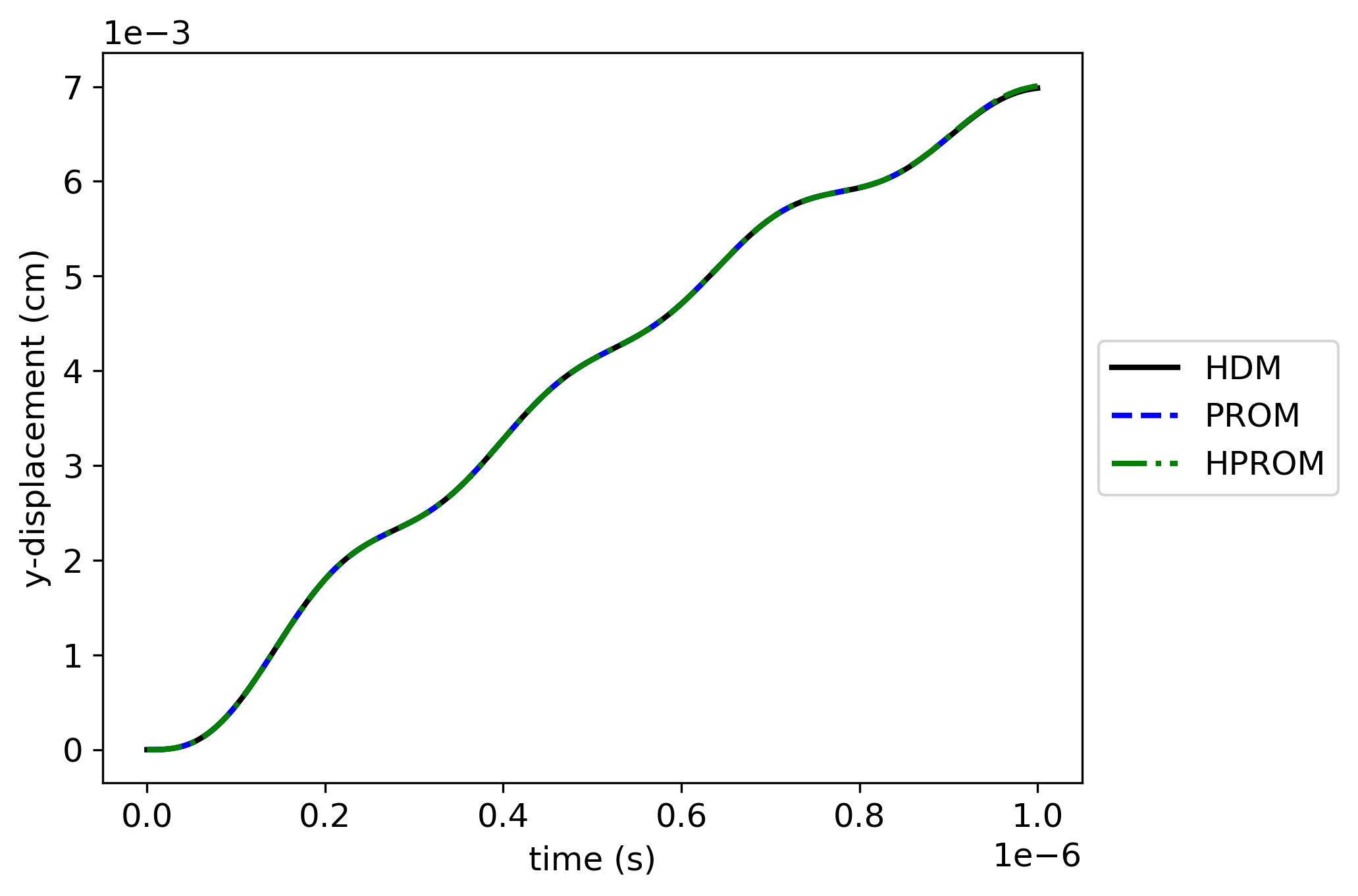}
\caption{Displacement}
\label{fig:p2_y_disp}
\end{subfigure}
\begin{subfigure}[b]{0.48\textwidth}
\centering
\includegraphics[width=\textwidth]{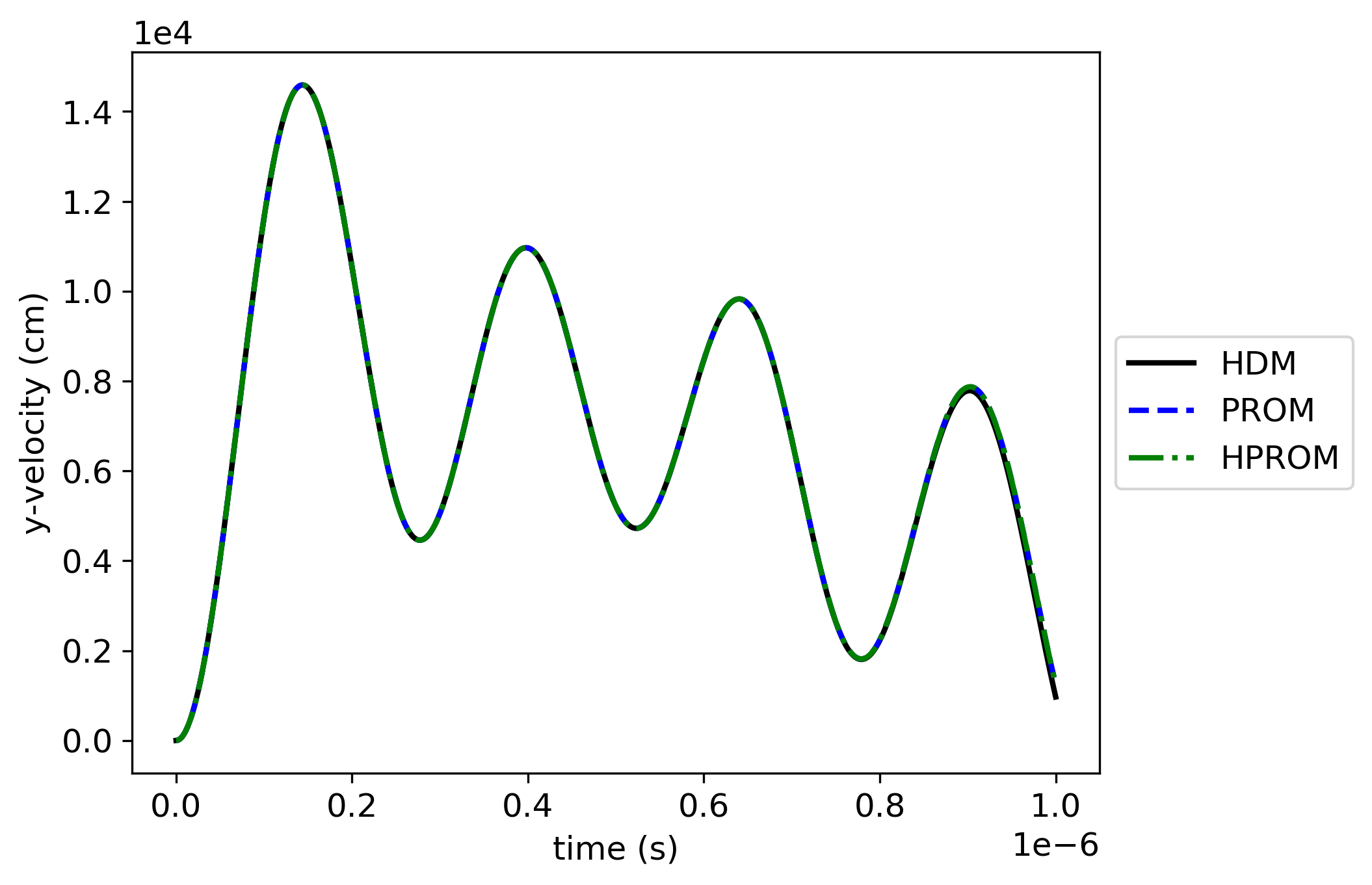}
\caption{Velocity}
\label{fig:p2_y_velo}
\end{subfigure}
\caption{Impact on a rigid wall of a nonlinear elastic cylinder with a composite microstructure -- time-histories of
         a displacement DOF and a velocity DOF at a node of the mesh underlying the macroscale HDM computed
         using: the macroscale and microscale HDMs; and the in-situ, adaptive, microscale framework for
         PMOR with and without hyperreduction.}
\label{fig:p2_history}
\end{figure}

Figure \ref{fig:p2_history} displays two sample time-histories computed for a displacement DOF and a
velocity DOF at a same node of the macroscale FE model. It shows that both adaptive PROM and HPROM reproduce
extremely accurately the predictions obtained using the macroscale and microscale HDMs. This impressive
accuracy of the proposed adaptive PMOR framework is confirmed by the quantitative error results reported
in Table \ref{tab:p2}. This table also shows that for the application considered herein, the in-situ,
adaptive, microscale framework for PMOR with and without hyperreduction accelerates the time to solution
of the simulation based on the macroscale and microscale HDMs by a factor greater than 7. Again, this
is an excellent performance result considering that it accounts for the computational time elapsed in
the processing of the macroscale model, which is neither reduced nor hyperreduced, and accounts for all
computational overhead associated with the PMOR process -- including adaptation.

\begin{table}[h!]
\centering
\caption{Impact on a rigid wall of a nonlinear elastic cylinder with a composite microstructure: performance
         assessment of the in-situ, adaptive, microscale framework for PMOR with and without hyperreduction.}
\begin{tabular}{ccc|ccc}
\hiderowcolors
\hline 
Model & \ Wall-clock time (s) & Speedup factor & Error measure $e$ & PROM & HPROM\\
\hline 
HDM & $2.44 \times 10^4$ & N.A. & $x$-displacement (\%) & $0.00185$ & $0.000467$\\
PROM & $1.60 \times 10^4$ & $1.52$ & $y$-displacement (\%) & $0.000165$ & $0.000350$\\
HPROM & $3.34 \times 10^3$ & $7.3$ & $x$-velocity (\%) & $0.000449$ & $0.00851$\\
 & & & y-velocity (\%) & $0.000430$ & $0.00991$\\
\hline
\end{tabular}
\label{tab:p2}
\end{table}

\newpage

\section{Conclusions}

An in-situ, adaptive, microscale Projection-based Model Order Reduction (PMOR) framework is proposed
in this paper for accelerating the simulation of the nonlinear dynamics of multiscale structural and
solid mechanics systems. It is designed to avoid by construction the weakness of extrapolation typically
found in nonadaptive training methods. This framework is based on the parameterization of the deformation
gradient for the purpose of PMOR at all but the coarsest scale, and on the concept of a database of local
Reduced-Order Bases (ROBs) where locality is measured in the deformation gradient parameter space. It
achieves accuracy by sampling new high-dimensional information as needed, updating on-the-fly a constructed
ROB, and approximating the solution along its trajectory using a handbook of most-appropriate local ROBs.
It achieves computational efficiency by splitting an updated ROB that becomes larger than desired into
two local ones, each capable of capturing only the local behavior of a nonlinear mesostructure or microstructure,
and incorporating hyperreduction (albeit without adaptivity). If the macroscale computational model is
not a parametric one, it is not reduced. Even in this case -- that is, when the macroscale computational
model is high-dimensional and its computational cost is accounted for -- and when all computational overhead
associated with performing PMOR and its adaptation are also accounted for, the proposed in-situ, adaptive,
PMOR framework accelerates three-dimensional, nonlinear, dynamic, multiscale simulations of the academic
type by roughly an order of magnitude without compromising accuracy. Larger speedup factors can be expected
for more realistic, large-scale problems and when adaptivity is extended to hyperreduction.

While described in the context of nonlinear, multiscale, dynamic problems in solid mechanics and structural
dynamics, the proposed in-situ, adaptive, PMOR framework is sufficiently comprehensive to be tailorable
to many other applications. Indeed, conventional offline training of ROBs in a predetermined region of
a parameter space typically leads to parametric reduced-order models that are bound to perform some form
of extrapolation, whenever they are exercised in a region of the parameter space that was not explored
during training. The proposed framework is designed specifically to remedy this issue.

\section*{acknowledgments}

The authors acknowledge partial support by the Air Force Office of Scientific Research under grant FA9550-17-1-0182,
and partial support by a research grant from the King Abdulaziz City for Science and Technology (KACST).
This document does not necessarily reflect the position of these institutions, and no official endorsement
should be inferred.

\bibliography{reference}

\begin{thebibliography}{25}
\providecommand{\natexlab}[1]{#1}
\providecommand{\url}[1]{\texttt{#1}}
\providecommand{\urlprefix}{}

\bibitem[{Qiao et~al.(2008)Qiao, Pizhong and Yang, Mijia and Bobaru,
  Florin}]{qiao2008impact}
Qiao P, Yang M, Bobaru F.
\newblock Impact mechanics and high-energy absorbing materials.
\newblock Journal of Aerospace Engineering 2008;21(4):235--248.

\bibitem[{Ma and Zhu(2017)Ma, Evan and Zhu, Ting}]{ma2017towards}
Ma E, Zhu T.
\newblock Towards strength--ductility synergy through the design of
  heterogeneous nanostructures in metals.
\newblock Materials Today 2017;20(6):323--331.

\bibitem[{Hoffler et~al.(2000)Hoffler, CE and Moore, KE and Kozloff, K and
  Zysset, PK and Brown, MB and Goldstein, SA}]{hoffler2000heterogeneity}
Hoffler C, Moore K, Kozloff K, Zysset P, Brown M, Goldstein S.
\newblock Heterogeneity of bone lamellar-level elastic moduli.
\newblock Bone 2000;26(6):603--609.

\bibitem[{Mei and Vernescu(2010)Mei, Chiang C and Vernescu,
  Bogdan}]{mei2010homogenization}
Mei CC, Vernescu B.
\newblock Homogenization methods for multiscale mechanics.
\newblock World scientific; 2010.

\bibitem[{Smit et~al.(1998)Smit, RJM and Brekelmans, WAM and Meijer,
  HEH}]{smit1998prediction}
Smit R, Brekelmans W, Meijer H.
\newblock Prediction of the mechanical behavior of nonlinear heterogeneous
  systems by multi-level finite element modeling.
\newblock Computer methods in applied mechanics and engineering
  1998;155(1-2):181--192.

\bibitem[{Miehe et~al.(1999)Miehe, Christian and Schr{\"o}der, J{\"o}rg and
  Schotte, Jan}]{miehe1999computational}
Miehe C, Schr{\"o}der J, Schotte J.
\newblock Computational homogenization analysis in finite plasticity simulation
  of texture development in polycrystalline materials.
\newblock Computer methods in applied mechanics and engineering
  1999;171(3-4):387--418.

\bibitem[{Feyel and Chaboche(2000)Feyel, Fr{\'e}d{\'e}ric and Chaboche,
  Jean-Louis}]{feyel2000fe2}
Feyel F, Chaboche JL.
\newblock FE2 multiscale approach for modelling the elastoviscoplastic
  behaviour of long fibre SiC/Ti composite materials.
\newblock Computer methods in applied mechanics and engineering
  2000;183(3-4):309--330.

\bibitem[{Kouznetsova et~al.(2001)Kouznetsova, V and Brekelmans, WAM and
  Baaijens, FPT}]{kouznetsova2001approach}
Kouznetsova V, Brekelmans W, Baaijens F.
\newblock An approach to micro-macro modeling of heterogeneous materials.
\newblock Computational mechanics 2001;27(1):37--48.

\bibitem[{Zahr et~al.(2017)Zahr, Matthew J and Avery, Philip and Farhat,
  Charbel}]{zahr2017multilevel}
Zahr MJ, Avery P, Farhat C.
\newblock A multilevel projection-based model order reduction framework for
  nonlinear dynamic multiscale problems in structural and solid mechanics.
\newblock International Journal for Numerical Methods in Engineering
  2017;112(8):855--881.

\bibitem[{Sirovich(1987)Sirovich, Lawrence}]{sirovich1987turbulence}
Sirovich L.
\newblock Turbulence and the dynamics of coherent structures. I. Coherent
  structures.
\newblock Quarterly of applied mathematics 1987;45(3):561--571.

\bibitem[{LeGresley and Alonso(2000)LeGresley, Patrick and Alonso,
  Juan}]{legresley2000airfoil}
LeGresley P, Alonso J.
\newblock Airfoil design optimization using reduced order models based on
  proper orthogonal decomposition.
\newblock In: Fluids 2000 conference and exhibit; 2000. p. 2545.

\bibitem[{Manzoni et~al.(2012)Manzoni, Andrea and Quarteroni, Alfio and Rozza,
  Gianluigi}]{manzoni2012shape}
Manzoni A, Quarteroni A, Rozza G.
\newblock Shape optimization for viscous flows by reduced basis methods and
  free-form deformation.
\newblock International Journal for Numerical Methods in Fluids
  2012;70(5):646--670.

\bibitem[{Farhat et~al.(2018{\natexlab{a}})Farhat, Charbel and Bos, Adrien and
  Tezaur, Radek and Chapman, Todd and Avery, Philip and Soize,
  Christian}]{farhat2018stochastic}
Farhat C, Bos A, Tezaur R, Chapman T, Avery P, Soize C.
\newblock A stochastic projection-based hyperreduced order model for model-form
  uncertainties in vibration analysis.
\newblock In: 2018 AIAA Non-Deterministic Approaches Conference; 2018. p. 1410.

\bibitem[{Farhat et~al.(2018{\natexlab{b}})Farhat, Charbel and Bos, Adrien and
  Avery, Philip and Soize, Christian}]{farhat2018modeling}
Farhat C, Bos A, Avery P, Soize C.
\newblock Modeling and quantification of model-form uncertainties in eigenvalue
  computations using a stochastic reduced model.
\newblock AIAA Journal 2018;56(3):1198--1210.

\bibitem[{Hetmaniuk et~al.(2012)Hetmaniuk, U and Tezaur, R and Farhat,
  C}]{hetmaniuk2012review}
Hetmaniuk U, Tezaur R, Farhat C.
\newblock Review and assessment of interpolatory model order reduction methods
  for frequency response structural dynamics and acoustics problems.
\newblock International Journal for Numerical Methods in Engineering
  2012;90(13):1636--1662.

\bibitem[{Amsallem et~al.(2010)Amsallem, David and Cortial, Julien and Farhat,
  Charbel}]{amsallem2010towards}
Amsallem D, Cortial J, Farhat C.
\newblock Towards real-time computational-fluid-dynamics-based aeroelastic
  computations using a database of reduced-order information.
\newblock AIAA journal 2010;48(9):2029--2037.

\bibitem[{Yvonnet and He(2007)Yvonnet, Julien and He, Q-C}]{yvonnet2007reduced}
Yvonnet J, He QC.
\newblock The reduced model multiscale method (R3M) for the non-linear
  homogenization of hyperelastic media at finite strains.
\newblock Journal of Computational Physics 2007;223(1):341--368.

\bibitem[{Monteiro et~al.(2008)Monteiro, E and Yvonnet, Julien and He,
  Qi-Chang}]{monteiro2008computational}
Monteiro E, Yvonnet J, He QC.
\newblock Computational homogenization for nonlinear conduction in
  heterogeneous materials using model reduction.
\newblock Computational Materials Science 2008;42(4):704--712.

\bibitem[{Farhat et~al.(2014)Farhat, Charbel and Avery, Philip and Chapman,
  Todd and Cortial, Julien}]{farhat2014dimensional}
Farhat C, Avery P, Chapman T, Cortial J.
\newblock Dimensional reduction of nonlinear finite element dynamic models with
  finite rotations and energy-based mesh sampling and weighting for
  computational efficiency.
\newblock International Journal for Numerical Methods in Engineering
  2014;98(9):625--662.

\bibitem[{Farhat et~al.(2015)Farhat, Charbel and Chapman, Todd and Avery,
  Philip}]{farhat2015structure}
Farhat C, Chapman T, Avery P.
\newblock Structure-preserving, stability, and accuracy properties of the
  energy-conserving sampling and weighting method for the hyper reduction of
  nonlinear finite element dynamic models.
\newblock International Journal for Numerical Methods in Engineering
  2015;102(5):1077--1110.

\bibitem[{Amsallem et~al.(2012)Amsallem, David and Zahr, Matthew J and Farhat,
  Charbel}]{amsallem2012nonlinear}
Amsallem D, Zahr MJ, Farhat C.
\newblock Nonlinear model order reduction based on local reduced-order bases.
\newblock International Journal for Numerical Methods in Engineering
  2012;92(10):891--916.

\bibitem[{Lawson and Hanson(1995)Lawson, Charles L and Hanson, Richard
  J}]{lawson1995solving}
Lawson CL, Hanson RJ.
\newblock Solving least squares problems, vol.~15.
\newblock Siam; 1995.

\bibitem[{Chapman et~al.(2017)Chapman, Todd and Avery, Philip and Collins, Pat
  and Farhat, Charbel}]{chapman2017accelerated}
Chapman T, Avery P, Collins P, Farhat C.
\newblock Accelerated mesh sampling for the hyper reduction of nonlinear
  computational models.
\newblock International Journal for Numerical Methods in Engineering
  2017;109(12):1623--1654.

\bibitem[{Zahr and Farhat(2015)Zahr, Matthew J and Farhat,
  Charbel}]{zahr2015progressive}
Zahr MJ, Farhat C.
\newblock Progressive construction of a parametric reduced-order model for
  PDE-constrained optimization.
\newblock International Journal for Numerical Methods in Engineering
  2015;102(5):1111--1135.

\bibitem[{Taylor(1948)Taylor, Geoffrey Ingram}]{taylor1948use}
Taylor GI.
\newblock The use of flat-ended projectiles for determining dynamic yield
  stress I. Theoretical considerations.
\newblock Proceedings of the Royal Society of London Series A Mathematical and
  Physical Sciences 1948;194(1038):289--299.

\end{thebibliography}
\end{document}